\documentclass[a4paper,10pt,twoside]{article}
\usepackage{amsmath,amssymb,amsthm,eucal,graphics}
\usepackage[cp1251]{inputenc}
\usepackage[russian]{babel}
\usepackage{graphicx}
\usepackage{amsfonts,amssymb,eucal}
\usepackage{amsbsy}
\usepackage{latexsym}
\usepackage{tabularx}

\makeatletter \gdef\firstpage{1}

\def\frsthdr{ }

\def\firstpageone{0\thepage}
\def\firstpagetwo{00\thepage}
\def\firstpagethree{000\thepage}
\def\firstpagemark{\ifnum\firstpage <10  \firstpageone
\else\ifnum\firstpage<100 \firstpagetwo \else \ifnum\firstpage
<1000 \firstpagethree \else \firstpageone\fi\fi\fi}

\def\footline{\ifnum\thepage=\firstpage \footlineone
\else\footlinetwo\fi}
\def\footlineone{ }
\def\footlinetwo{}

\def\titles{Finite almost simple groups}
\def\authors{ I.\,B.\,Gorskov, N.\,V.\,Maslova}

\def\oddhedr{\ifnum\thepage=\firstpage \firsthdr \else \odhdr \fi}
\def\firsthdr{\hspace{\fill} \sl \frsthdr \hspace{\fill}\hbox{}}
\def\odhdr{\hspace{\fill}\sl\rightmark \titles \hspace{\fill}
\rm \thepage}

\def\evnhedr{\ifnum\thepage=\firstpage \firsthdr \else \evhdr \fi}
\def\evhdr{\noindent \rm \thepage\hspace*{\fill} \sl\leftmark
\authors \hspace*{\fill}\hbox{}}

\def\ps@newpstyle{\def\@oddhead{
\hspace{-0.65em} \vbox{\oddhedr\vskip 1mm \hrule width \textwidth}
}
\def\@evenhead{
\hspace{-0.65em} \vbox{\evnhedr\vskip 1mm \hrule width \textwidth}
}\textsc{}
\def\@oddfoot{\footline}
\def\@evenfoot{\@oddfoot}
}

\def\refer
{
\begin{small}
\begin{enumerate}
\leftmargin=0pt \rightmargin=0pt \itemsep=0pt
\parsep=0pt
}

\def\endref{\end{enumerate}\end{small} }

\begin{document}

\setcounter{page}{\firstpage} 

\Russian \sloppy \rm


\newcommand{\NN}{\mathbb{N}}
\newcommand{\ZZ}{\mathbb{Z}}
\newcommand{\RR}{\mathbb{R}}

\newcommand{\inN}{\in\NN}
\newcommand{\inZ}{\in\ZZ}
\newcommand{\inR}{\in\RR}

\newcommand{\alf}{\alpha}
\newcommand{\bet}{\beta}
\newcommand{\eps}{\varepsilon}
\newcommand{\lam}{\lambda}
\newcommand{\ups}{\upsilon}
\newcommand{\sgn}{\rm sign}

\newtheorem{teo}{\bf ╥┼╬╨┼╠└}
\newtheorem{lem}{\bf ╦┼╠╠└}
\newtheorem{hyp}{\bf ╧╨┼─╧╬╦╬╞┼═╚┼}
\newtheorem{prp}{\bf ╧╨┼─╦╬╞┼═╚┼}
\newtheorem{ass}{\bf ╙╥┬┼╨╞─┼═╚┼}
\newtheorem{cor}{\bf ╤╦┼─╤╥┬╚┼}

\newcommand{\pr}{\par\mbox{─╬╩└╟└╥┼╦▄╤╥┬╬.}\ \ }

\renewcommand{\theteo}{\arabic{teo}}
\renewcommand{\theprp}{\arabic{prp}}
\renewcommand{\theass}{\arabic{ass}}
\renewcommand{\thelem}{\arabic{lem}}
\renewcommand{\thecor}{\arabic{cor}}
\renewcommand{\thehyp}{\arabic{hyp}}

\begin{center} FINITE ALMOST SIMPLE GROUPS WHOSE GRUENBERG--KEGEL GRAPHS \\ COINCIDE WITH GRUENBERG--KEGEL GRAPHS OF SOLVABLE GROUPS\\
(in Russian, English abstract)\\
I.\,B.\,Gorshkov$^{1}$, N.\,V.\,Maslova$^{1,2}$\\
${^1}$Krasovskii Institute of Mathematics and Mechanics UB RAS, \\
${^2}$Ural Federal University\\
ilygor8@gmail.com, butterson@mail.ru
\end{center}

\allowdisplaybreaks
\renewcommand{\tablename}{}

\footnote{\small { The work is supported by the grant of the President of the Russian Federation for young scientists (grant no. MK-6118.2016.1), by the Integrated Program for Fundamental Research of the Ural Branch of the Russian Academy of Sciences (project no. 15-16-1-5), by the Program for State Support of Leading Universities of the Russian Federation (agreement
no. 02.A03.21.0006 of August 27, 2013). The first author is supported by the grant FAPESP-2014/08964-1. The second author is a winner of the competition of young mathematicians of the Dmitry Zimin Foundation ``Dynasty'' (in 2013 year).}}

{\small \small \emph{Abstract.} Let~$G$ be a finite group. Denote by $\pi(G)$ the set of all prime divisors of the
order of $G$ and by $\omega (G)$
the {\it spectrum} of~$G$, i.e. the set of all its element orders. The
set $\omega(G)$ defines the {\it Gruenberg--Kegel graph} (or the {\it prime graph})
$\Gamma(G)$ of~$G$; in this graph the vertex set is $\pi(G)$ and different vertices~$p$ and~$q$ are adjacent if and only if
$pq\in\omega (G)$.

In [J. Group Theory. Vol.~2, no.~2 (1999). P. 157--172] Lucido proved that if $G$ is a finite solvable group then $\Gamma(G)$ does not contain a $3$-coclique. In [J. Algebra. Vol.~442 (2015). P.~397--422] Gruber et. al. proved that a graph is isomorphic to the Gruenberg--Kegel graph of a finite solvable group if and only if its complement is $3$-colorable and triangle free.
In [Proc. Steklov Inst. Math. Vol.~283, suppl. 1 (2013). P.~139Ц145] Zinov'eva and Mazurov described finite simple groups whose Gruenberg--Kegel graphs coincide with Gruenberg--Kegel graphs of Frobenius groups or $2$-Frobenius groups. It is not difficult to prove these simple groups are exactly all finite simple groups whose Gruenberg--Kegel graphs coincide with Gruenberg--Kegel graphs of finite solvable groups. Our aim is to obtain similar results for finite almost simple groups. In this paper, we prove the following theorem.

\smallskip

{\bf Theorem 1.} Let $G$ be a finite almost simple group. Then the following conditions are equivalent{\rm:}

$(1)$ $\Gamma(G)$ does not contain a $3$-coclique{\rm;}

$(2)$ $\Gamma(G)$ is isomorphic to the Gruenberg--Kegel graph of a finite solvable group{\rm;}

$(3)$ $\Gamma(G)$ is equal to the Gruenberg--Kegel graph of an appropriate finite solvable group.

\smallskip


The following question naturally arises: is there a graph $\Gamma$ without $3$-cocliques such that the complement of $\Gamma$ is not $3$-colorable and $\Gamma$ is isomorphic to the Gruenberg--Kegel graph of an appropriate finite non-solvable group?

\smallskip

Moreover, we obtain an explicit description of finite almost simple groups whose Gruenberg--Kegel graphs coincide with Gruenberg--Kegel graphs of finite solvable groups. We prove the following theorem.

\smallskip

{\bf Theorem 2.} Let $G$ be a finite almost simple group such that $S=Soc(G)$ is isomorphic to one of the following simple groups: $A_n$, where $n\geq 5$; $PSL_n(q)$, where $n \ge 2$ and $(n,q) \not = (2,2),(2,3)$; $PSU_n(q)$, where $n \ge 3$ and $(n,q) \not = (3,2)$; $PSp_n(q)$, where $n \ge 4$ is even; $P\Omega_n(q)$, where $n \ge 7$ is odd; $P\Omega_n^\pm(q)$, where $n \ge 8$ is even; an exceptional group of Lie type over the field of order $q$; a sporadic simple group. Let $f$ be the standard field automorphism of $S$ and $g$ be the standard graph automorphism of $S$. Then $\Gamma(G)$ does not contain a $3$-coclique if and only if $\pi(G)=\pi(S)$ and one of the following conditions holds{\rm:}

$(1)$ $S$ is isomorphic to one of the following groups{\rm:} $A_9$, $A_{10}$, $A_{12}$, $PSU_3(9)$, $PSU_4(2)$, $PSp_6(2)$, $P\Omega_8^+(2)$, ${^3}D_4(2)${\rm;}

$(2)$ $G$ is isomorphic to one of the following groups{\rm:} $S_5$, $S_6$, $PGL_2(9)$, $M_{10}$, $Aut(A_6)$, $S_8$, $Aut(PSL_2(8))$,
$Aut(PSL_3(2))$, $PGL_3(4)\langle f \rangle$, $PGL_3(4)\langle g \rangle$, $Aut(PSL_3(4))$, $PSL_4(4)\langle f \rangle$, $PSL_4(4)\langle g \rangle$, $Aut(PSL_4(4))$, $Aut(PSU_5(2))${\rm;}

$(3)$ $G \cong PGL_2(p)$, where $p$ is either a Fermat prime or a Mersenne prime{\rm;}

$(4)$ $S \cong PSL_2(2^m)$, where $m \ge 4$ is even and $\{2\} \subseteq \pi(G/S)${\rm;}

$(5)$ $S \cong PSL_3(p)$, where $p$ is a Mersenne prime and $(p-1)_3\not =3${\rm;}

$(6)$ $S \cong PSL_3(p)$, where $p$ is either a Fermat prime or a Mersenne prime, $(p-1)_3=3$ and $Inndiag(S) \le G \le Aut(S)${\rm;}

$(7)$ $S \cong PSL_3(2^m)$, where $m \ge 3$, $(2^m-1)_3=3$ and $Inndiag(S)\langle g\rangle \le G \le Aut(PSL_3(2^m))${\rm;}

$(8)$ $S \cong PSL_3(2^m)$, where $m \ge 3$, $(2^m-1)_3\not =3$ and $S\langle g\rangle \le G \le Aut(PSL_3(2^m))${\rm;}

$(9)$ $S \cong PSL_4(2^m)$, where $m \ge 3$ and $S\langle g\rangle \le G \le Aut(PSL_4(2^m))${\rm;}

$(10)$ $S \cong PSU_3(p)$, where $p$ is a Fermat prime and $(p+1)_3\not =3${\rm;}

$(11)$ $S \cong PSU_3(p)$, where $p$ is a Fermat prime, $(p+1)_3=3$ and $Inndiag(S) \le G \le Aut(S)${\rm;}

$(12)$ $S \cong PSU_3(2^m)$, where $m \ge 2$, $(2^m-1)_3=3$, $\{2\} \subseteq \pi(G/S)$ and $Inndiag(S) \le G \le Aut(S)${\rm;}

$(13)$ $S \cong PSU_3(2^m)$, where $m \ge 2$, $(2^m-1)_3\not =3$ and $\{2\} \subseteq \pi(G/S)${\rm;}

$(14)$ $S \cong PSU_4(2^m)$, where $m \ge 2$ and $\{2\} \subseteq \pi(G/S)${\rm;}

$(15)$ $S \cong PSp_4(q)$.

\smallskip

\emph{Keywords:} finite group, solvable group, almost simple group, Gruenberg--Kegel graph (prime graph).}

\section*{┬тхфхэшх}

╠√ сєфхь єяюЄЁхсы Є№ ЄхЁьшэ <<уЁєяяр>> т чэрўхэшш <<ъюэхўэр  уЁєяяр>> ш ЄхЁьшэ <<уЁрЇ>> т чэрўхэшш <<ъюэхўэ√щ уЁрЇ схч яхЄхы№ ш ъЁрЄэ√ї ЁхсхЁ>>.

{\it ╤яхъЄЁюь} уЁєяя√ $G$ эрч√трхЄё  ьэюцхёЄтю $\omega(G)$ яюЁ фъют тёхї хх
¤ыхьхэЄют. ╠эюцхёЄтю тёхї яЁюёЄ√ї ўшёхы, тїюф ∙шї т ёяхъЄЁ уЁєяя√ $G$, эрч√трхЄё  {\it яЁюёЄ√ь ёяхъЄЁюь уЁєяя√ $G$} ш юсючэрўрхЄё  ўхЁхч $\pi(G)$. ╤яхъЄЁ $\omega(G)$ юяЁхфхы хЄ {\it уЁрЇ ├Ё■эсхЁур--╩хухы } (шыш {\it уЁрЇ яЁюёЄ√ї ўшёхы}) $\Gamma(G)$  уЁєяя√ $G$, тхЁ°шэрьш ъюЄюЁюую  ты ■Єё  яЁюёЄ√х ўшёыр шч $\pi(G)$, ш фтх Ёрчышўэ√х тхЁ°шэ√ $r$ ш $s$ ёьхцэ√ Єюуфр ш Єюы№ъю Єюуфр, ъюуфр ўшёыю $rs$ яЁшэрфыхцшЄ ьэюцхёЄтє $\omega(G)$. ╬сючэрўшь ўхЁхч $t(G)$ эршсюы№°хх ўшёыю яЁюёЄ√ї ўшёхы т $\pi(G)$, яюярЁэю
эх ёьхцэ√ї т $\Gamma(G)$.

═ряюьэшь, ўЄю яюф ЁртхэёЄтюь уЁрЇют ├Ё■эсхЁур--╩хухы  уЁєяя ь√ яюэшьрхь ёютярфхэшх шї ьэюцхёЄт тхЁ°шэ (Єю хёЄ№ яЁюёЄ√ї ёяхъЄЁют уЁєяя) ш ьэюцхёЄт ЁхсхЁ. ╧юф шчюьюЁЇшчьюь уЁрЇют ├Ё■эсхЁур--╩хухы  уЁєяя ь√ яюэшьрхь шї шчюьюЁЇшчь ъръ рсёЄЁръЄэ√ї уЁрЇют.

┬ \cite[ЄхюЁхьр~1]{Lucido2}*\footnote{*═рёЄю ∙шщ Ёхчєы№ЄрЄ с√ы яюыєўхэ ╠.╤. ╦єўшфю т 1999 у., юфэръю юэ эряЁ ьє■ ёыхфєхЄ шч сюыхх Ёрээшї Ёхчєы№ЄрЄют ├. ╒шуьрэр \cite[ЄхюЁхьр 1]{Hig1}  ш ╘. ╒юыыр \cite[ЄхюЁхьр 6.4.1]{Gorenstein}.} с√ыю чрьхўхэю, ўЄю уЁрЇ√ ├Ё■эсхЁур--╩хухы  ЁрчЁх°шь√ї уЁєяя эх ёюфхЁцрЄ $3$-ъюъышъ, шэ√ьш ёыютрьш $t(H) \le 2$ фы  ы■сющ ЁрчЁх°шьющ уЁєяя√ $H$. ┬ \cite{Keller} с√ыю яюърчрэю, ўЄю уЁрЇ шчюьюЁЇхэ уЁрЇє ├Ё■эсхЁур--╩хухы  ЁрчЁх°шьющ уЁєяя√ Єюуфр ш Єюы№ъю Єюуфр, ъюуфр юэ эх ёюфхЁцшЄ $3$-ъюъышъ ш хую фюяюыэхэшх $3$-ЁрёъЁр°штрхью. ┬ \cite[ЄхюЁхь√~1,~3]{ZinMaz} юяшёрэ√ яЁюёЄ√х уЁєяя√, уЁрЇ√ ├Ё■эсхЁур--╩хухы  ъюЄюЁ√ї Ёртэ√ уЁрЇрь ├Ё■эсхЁур--╩хухы  ЁрчЁх°шь√ї уЁєяя ╘Ёюсхэшєёр ш фтющэ√ї уЁєяя ╘Ёюсхэшєёр. ╚ч \cite[ЄрсышЎ√~2,3]{Maslova_VasVd2005} ш \cite[ЄрсышЎ√~2,~3,~4]{Maslova_VasVd2011} ёыхфєхЄ, ўЄю уЁєяярьш шч \cite[ЄхюЁхь√~1,~3]{ZinMaz} шёўхЁя√тр■Єё  тёх яЁюёЄ√х уЁєяя√, уЁрЇ√ ├Ё■эсхЁур--╩хухы  ъюЄюЁ√ї Ёртэ√ уЁрЇрь ├Ё■эсхЁур--╩хухы  ЁрчЁх°шь√ї уЁєяя. ╓хы№ эрёЄю ∙хщ ЁрсюЄ√ --- яюыєўшЄ№ яюфюсэюх юяшёрэшх фы  яюўЄш яЁюёЄ√ї уЁєяя. ═рьш фюърчрэр ёыхфє■∙р 

\begin{teo}\label{ASasS} ╧єёЄ№ $G$ --- яюўЄш яЁюёЄр  уЁєяяр. ╥юуфр ёыхфє■∙шх єёыютш  ¤ътштрыхэЄэ√{\rm:}

$(1)$ $\Gamma(G)$ эх ёюфхЁцшЄ $3$-ъюъышъ{\rm;}

$(2)$ $\Gamma(G)$ шчюьюЁЇхэ уЁрЇє ├Ё■эсхЁур--╩хухы  эхъюЄюЁющ ЁрчЁх°шьющ уЁєяя√{\rm;}

$(3)$ $\Gamma(G)$ Ёртхэ уЁрЇє ├Ё■эсхЁур--╩хухы  эхъюЄюЁющ ЁрчЁх°шьющ уЁєяя√.

\end{teo}

┼ёЄхёЄтхээ√ь юсЁрчюь тючэшърхЄ ёыхфє■∙шщ

\medskip

{\bf ┬╬╧╨╬╤.} {\it ╤є∙хёЄтєхЄ ыш уЁрЇ схч $3$-ъюъышъ, ъюЄюЁ√щ шчюьюЁЇхэ уЁрЇє ├Ё■эсхЁур--╩хухы  эхъюЄюЁющ
эхЁрчЁх°шьющ уЁєяя√, ш эх шчюьюЁЇхэ уЁрЇє ├Ё■эсхЁур--╩хухы  эшъръющ ЁрчЁх°шьющ уЁєяя√{\rm?}}
\medskip

═ряюьэшь, ўЄю ўхЁхч $Soc(G)$ юсючэрўрхЄё  Ўюъюы№ уЁєяя√ $G$ (яюфуЁєяяр, яюЁюцфхээр  тёхьш хх ьшэшьры№э√ьш эххфшэшўэ√ьш эюЁьры№э√ьш яюфуЁєяярьш). ╧єёЄ№ $q=p^m$ ш $G$ --- яюўЄш яЁюёЄр  уЁєяяр Єрър , ўЄю яюфуЁєяяр $S=Soc(G)$ шчюьюЁЇэр юфэющ шч яЁюёЄ√ї уЁєяя: $A_n$, уфх $n\geq 5$; $PSL_n(q)$, уфх $n \ge 2$ ш $(n,q) \not = (2,2),(2,3)$; $PSU_n(q)$, уфх $n \ge 3$ ш $(n,q) \not = (3,2)$; $PSp_n(q)$, уфх $n \ge 4$ ўхЄэю; $P\Omega_n(q)$, уфх $n \ge 7$ эхўхЄэю; $P\Omega_n^\pm(q)$, уфх $n \ge 8$ ўхЄэю; яЁюёЄющ шёъы■ўшЄхы№эющ уЁєяях ышхтр Єшяр эрф яюыхь яюЁ фър $q$ шыш юфэющ шч ъюэхўэ√ї яЁюёЄ√ї ёяюЁрфшўхёъшї уЁєяя. ═хюсїюфшь√х ётхфхэш  ю ёЄЁюхэшш уЁєяя ртЄюьюЁЇшчьют ъюэхўэ√ї яЁюёЄ√ї ышэхщэ√ї ш єэшЄрЁэ√ї уЁєяя яЁштхфхэ√ фрыхх т ыхььрї \ref{AutPSL} ш \ref{AutPSU}. ╩ы■ўхт√ь фы  фюърчрЄхы№ёЄтр ЄхюЁхь√ \ref{ASasS}  ты хЄё  ёыхфє■∙шщ яюыєўхээ√щ эрьш Ёхчєы№ЄрЄ, шьх■∙шщ эхчртшёшь√щ шэЄхЁхё.
\smallskip

\begin{teo}\label{AlmostSimple} ╧єёЄ№ $G$ --- яюўЄш яЁюёЄр  уЁєяяр ш $S=Soc(G)$. ╥юуфр $\Gamma(G)$ эх ёюфхЁцшЄ $3$-ъюъышъ Єюуфр ш Єюы№ъю Єюуфр, ъюуфр т√яюыэ хЄё  юфэю шч ёыхфє■∙шї єЄтхЁцфхэшщ{\rm:}

$(1)$ уЁєяяр $S$ шчюьюЁЇэр юфэющ шч ёыхфє■∙шї уЁєяя{\rm:} $A_9$, $A_{10}$, $A_{12}$, $PSU_3(9)$, $PSU_4(2)$, $PSp_6(2)$, $P\Omega_8^+(2)$, ${^3}D_4(2)${\rm;}

$(2)$ уЁєяяр $G$ шчюьюЁЇэр юфэющ шч ёыхфє■∙шї уЁєяя{\rm:} $S_5$, $S_6$, $PGL_2(9)$, $M_{10}$, $Aut(A_6)$, $S_8$, $Aut(PSL_2(8))$,
$Aut(PSL_3(2))$, $PGL_3(4)\langle f \rangle$, $PGL_3(4)\langle g \rangle$, $Aut(PSL_3(4))$, $PSL_4(4)\langle f \rangle$, $PSL_4(4)\langle g \rangle$, $Aut(PSL_4(4))$, $Aut(PSU_5(2))${\rm;}

$(3)$ $G \cong PGL_2(p)$, $p$ --- яЁюёЄюх ўшёыю ╘хЁьр шыш ╠хЁёхэр{\rm;}

$(4)$ $S \cong PSL_2(2^m)$, $m \ge 4$ ўхЄэю ш $\{2\} \subseteq \pi(G/S) \subseteq \pi(S)${\rm;}

$(5)$ $S \cong PSL_3(p)$, $p$ --- яЁюёЄюх ўшёыю ╠хЁёхэр ш $(p-1)_3\not =3${\rm;}

$(6)$ $S \cong PSL_3(p)$, $p$ --- яЁюёЄюх ўшёыю ╠хЁёхэр, $(p-1)_3=3$ ш $Inndiag(S) \le G \le Aut(S)${\rm;}

$(7)$ $S \cong PSL_3(2^m)$, $m \ge 3$, $(2^m-1)_3=3$, $\pi(G) = \pi(S)$ ш $Inndiag(S)\langle g\rangle \le G \le Aut(S)${\rm;}

$(8)$ $S \cong PSL_3(2^m)$, $m \ge 3$, $(2^m-1)_3\not =3$, $\pi(G) = \pi(S)$ ш $S\langle g\rangle \le G \le Aut(S)${\rm;}

$(9)$ $S \cong PSL_4(2^m)$, $m \ge 3$, $\pi(G) = \pi(S)$ ш $S\langle g\rangle \le G \le Aut(S)${\rm;}

$(10)$ $S \cong PSU_3(p)$, $p$ --- яЁюёЄюх ўшёыю ╘хЁьр ш $(p+1)_3\not =3${\rm;}

$(11)$ $S \cong PSU_3(p)$, $p$ --- яЁюёЄюх ўшёыю ╘хЁьр, $(p+1)_3=3$ ш $Inndiag(S) \le G \le Aut(S)${\rm;}

$(12)$ $S \cong PSU_3(2^m)$, $m \ge 2$, $(2^m-1)_3=3$, $\{2\} \subseteq \pi(G/S) \subseteq \pi(S)$ ш $Inndiag(S) \le G \le Aut(S)${\rm;}

$(13)$ $S \cong PSU_3(2^m)$, $m \ge 2$, $(2^m-1)_3\not =3$ ш $\{2\} \subseteq \pi(G/S) \subseteq \pi(S)${\rm;}

$(14)$ $S \cong PSU_4(2^m)$, $m \ge 2$ ш $\{2\} \subseteq \pi(G/S) \subseteq \pi(S)${\rm;}

$(15)$ $S \cong PSp_4(q)$ ш $\pi(G)=\pi(S)$.

\end{teo}

{\bf ╟└╠┼╫└═╚┼.} ┼ёыш $G$ --- яюўЄш яЁюёЄр  уЁєяяр Єрър , ўЄю $\Gamma(G)$ эх ёфхЁцшЄ $3$-ъюъышъ, Єю $\pi(G)=\pi(Soc(G))$.

\section*{┬ёяюьюурЄхы№э√х Ёхчєы№ЄрЄ√}

═р°ш юсючэрўхэш  ш ЄхЁьшэюыюуш  т юёэютэюь ёЄрэфрЁЄэ√, шї ьюцэю
эрщЄш т \cite{Gorenstein,Har,Atlas,Maslova_GorLySo}.

┼ёыш $q$ Ч эрЄєЁры№эюх ўшёыю, $r$ Ч эхў╕Єэюх яЁюёЄюх ўшёыю ш $(q, r) = 1$,
Єю $e(r, q)$ юсючэрўрхЄ ьєы№ЄшяышърЄштэ√щ яюЁ фюъ ўшёыр $q$ яю ьюфєы■ $r$,
Ёртэ√щ ьшэшьры№эюьє эрЄєЁры№эюьє ўшёыє $m$, єфютыхЄтюЁ ■∙хьє єёыютш■ $q^m \equiv 1 \pmod r$.
─ы  эхў╕Єэюую $q$ яюыюцшь $e(2, q) = 1$, хёыш $q \equiv 1 \pmod 4$, ш $e(2, q) = 2$ т яЁюЄштэюь ёыєўрх.

\begin{lem}\label{Maslova_l2zh} {\rm(ёыхфёЄтшх ЄхюЁхь√ ╞шуьюэфш, \cite{Maslova_Zs})}
╧єёЄ№ $q$ --- эрЄєЁры№эюх ўшёыю, сюы№°хх $1$. ─ы  ы■сюую эрЄєЁры№эюую ўшёыр $m$ ёє∙хёЄтєхЄ яЁюёЄюх ўшёыю $r$,
єфютыхЄтюЁ ■∙хх ЁртхэёЄтє $e(r, q) = m$, чр шёъы■ўхэшхь ёыхфє■∙шї ёыєўрхт: $q = 2$ ш $m = 1$; $q = 3$ ш $m = 1$; $q = 2$ ш $m = 6$.
\end{lem}

╧ЁюёЄюх ўшёыю $r$, єфютыхЄтюЁ ■∙хх єёыютш■ $e(r, q) = m$, эрч√трхЄё  яЁшьшЄштэ√ь яЁюёЄ√ь фхышЄхыхь ўшёыр $q^m - 1$. ╧ю ыхььх~\ref{Maslova_l2zh} Єръюх ўшёыю ёє∙хёЄтєхЄ, чр шёъы■ўхэшхь ёыєўрхт, єяюь эєЄ√ї т ыхььх. ─ы  фрээюую $q$ юсючэрўшь ўхЁхч $R_m(q)$ ьэюцхёЄтю тёхї яЁшьшЄштэ√ї яЁюёЄ√ї фхышЄхыхщ ўшёыр $q^m-1$. ╦хуъю яюэ Є№, ўЄю $R_i(q) \cap R_j(q)=\emptyset$ яЁш ы■с√ї Ёрчышўэ√ї $i$ ш $j$, ш хёыш $q$ --- ёЄхяхэ№ яЁюёЄюую ўшёыр, Єю $R_i(q) \cap \pi(q) = \emptyset$ фы  ы■сюую $i$.

\begin{lem}[{\rm ыхььр ├хЁюэю \cite{Maslova_Gerono}}]\label{Maslova_Gerono} ╧єёЄ№ $p,~q$ --- яЁюёЄ√х ўшёыр Єръшх, ўЄю $p^a-q^b=1$ фы  эхъюЄюЁ√ї эрЄєЁры№э√ї ўшёхы $a,~b$. ╥юуфр ярЁр $(p^a,~q^b)$ Ёртэр $(3^2,~2^3),~(p,~2^b)$ шыш $(2^a,~q)$.
\end{lem}



\begin{lem}\label{solvable} ╧єёЄ№ $\pi$ --- ъюэхўэюх ьэюцхёЄтю яЁюёЄ√ї ўшёхы ш $\Gamma$ --- уЁрЇ, тхЁ°шэрьш ъюЄюЁюую  ты ■Єё  яЁюёЄ√х ўшёыр шч $\pi$. ┼ёыш $\pi =\pi_1 \cup \pi_2$, яЁшўхь $\pi_1$ ш $\pi_2$ --- эхяхЁхёхър■∙шхё  ъышъш т $\Gamma$, Єю ёє∙хёЄтєхЄ ЁрчЁх°шьр  уЁєяяр $H$ Єрър , ўЄю $\Gamma=\Gamma(H)$.
\end{lem}

\proof ╧єёЄ№ $\pi_1 = \{p_1, \ldots, p_m\}$, $\pi_2=\{q_1, \ldots, q_n\}$ ш $C = C_{q_1} \times \ldots \times C_{q_n}$ --- Ўшъышўхёър  уЁєяяр. ╧єёЄ№, схч юуЁрэшўхэш  юс∙эюёЄш, ўшёыю $p_1$ эх ёьхцэю т $\Gamma$ ё ўшёырьш шч ьэюцхёЄтр $\{q_1, \ldots q_k\}$ ш ёьхцэю ё ўшёырьш шч ьэюцхёЄтр $\{q_{k+1}, \ldots, q_n\}$. ┬тшфє ьрыющ ЄхюЁхь√ ╘хЁьр ёє∙хёЄтєхЄ эрЄєЁры№эюх ўшёыю $m_1$ Єръюх, ўЄю $p_1^{m_1}-1$ фхышЄё  эр яЁюшчтхфхэшх $q_1\cdot \ldots \cdot q_k$. ╥юуфр ьєы№ЄшяышърЄштэр  уЁєяяр $F_1^*$ яюы  $F_1$ яюЁ фър $p_1^{m_1}$ ёюфхЁцшЄ ¤ыхьхэЄ $x_1$ яюЁ фър $q_1\cdot \ldots \cdot q_k$. ╬яЁхфхышь фхщёЄтшх уЁєяя√ $C$ эр рффшЄштэющ уЁєяях $F_1^{+}$ яюы  $F_1$, яюырур , ўЄю фхщёЄтшх яюЁюцфр■∙хую ¤ыхьхэЄр яюфуЁєяя√ $C_{q_1} \times \ldots \times C_{q_k}$ эр $F_1^+$ --- ¤Єю єьэюцхэшх т яюых $F_1$ эр ¤ыхьхэЄ $x_1$, р фхщёЄтшх яюЁюцфр■∙хую ¤ыхьхэЄр яюфуЁєяя√ $C_{q_{k+1}} \times \ldots \times C_{q_n}$ эр $F_1^+$ ЄЁштшры№эю. └эрыюушўэ√ь юсЁрчюь юяЁхфхышь яюых $F_i$ ш фхщёЄтшх уЁєяя√ $C$ эр рффшЄштэющ уЁєяях $F_i^{+}$ яюы  $F_i$ фы  ърцфюую $p_i \in \pi_1$. ╧єёЄ№ уЁєяяр $H$ --- хёЄхёЄтхээюх яюыєяЁ ьюх яЁюшчтхфхэшх уЁєяя√ $F = F_1^+ \times \ldots F_m^+$ эр уЁєяяє $C$. ╥юуфр уЁєяяр $H$ ЁрчЁх°шьр ш уЁрЇ ├Ё■эсхЁур--╩хухы  $\Gamma(H)$ т ЄюўэюёЄш Ёртхэ уЁрЇє $\Gamma$. ╦хььр фюърчрэр.

\smallskip

╥ръцх эрь яюэрфюсшЄё  ёыхфє■∙р  ыхуъю фюърч√трхьр 
\begin{lem}\label{Maslova_nonadj} {\rm(ёь., эряЁшьхЁ, \cite[ыхььр~2]{Maslova_PrimeGM})} ╧єёЄ№ $L$ ---  уЁєяяр, $K \lhd L$, $r, s \in \pi(K)\setminus \pi(|L:K|)$. ╥юуфр хёыш яЁюёЄ√х ўшёыр $r$ ш $s$ эх ёьхцэ√ т $\Gamma(K)$, Єю юэш эх ёьхцэ√ т $\Gamma(L)$.
\end{lem}

\begin{lem}{\rm \cite[ЄхюЁхьр 4.5.1, яЁхфыюцхэш  2.5.12,\,4.9.1,\,4.9.2]{Maslova_GorLySo}}\label{AutPSL}
    ╧єёЄ№ $S = PSL_n(q)$, уфх $n\geq 2$, $q = p^m$ ш $(n,q) \not = (2,2),(2,3)$, $x$ --- ¤ыхьхэЄ яЁюёЄюую яюЁ фър $r$
    т $Aut(S)\setminus Inndiag(S)$ ш $S_x = O^{p'}(C_S(x))$. ╥юуфр ёяЁртхфышт√ ёыхфє■∙шх єЄтхЁцфхэш :

    $(1)$ $Aut(S)=Inndiag(S)\leftthreetimes (\Phi\times\langle g \rangle)$, уфх $Outdiag(S)\cong Z_{(n,\,q - 1)}$,
    $\Phi=\langle f \rangle\cong Aut(GF(q))\cong Z_m$ --- уЁєяяр яюыхт√ї ртЄюьюЁЇшчьют уЁєяя√ $S$, $g=1$ фы  $n=2$
    ш $g$ --- уЁрЇют√щ ртЄюьюЁЇшчь яюЁ фър $2$ уЁєяя√ $S$ фы  $n\geq 3$, $\Phi$ фхщёЄтєхЄ эр $Outdiag(S)$ Єръ цх, ъръ $Aut(GF(q))$
    фхщёЄтєхЄ эр ьєы№ЄшяышърЄштэє■ яюфуЁєяяє шч $GF(q)$ Єюую цх яюЁ фър, ўЄю ш $Outdiag(S)$, ш $g$ шэтхЁЄшЁєхЄ уЁєяяє $Outdiag(S)$;

    $(2)$ ┼ёыш $x\in Inndiag(S)\Phi$, Єю $r$ фхышЄ $m$, $S_x\cong PSL_n(q^{\frac{1}{r}})$ ш $C_{Inndiag(S)}(x)\cong Inndiag(S_x)$;

    $(3)$ ┼ёыш $n\geq 3$, $2$ фхышЄ $m$ ш $x\in Inndiag(S)f^{\frac{m}{2}}g$, Єю $r=2$, $S_x\cong PSU_n(q^{\frac{1}{2}})$, ш $C_{Inndiag(S)}(x)\cong Inndiag(S_x)$;

    $(4)$ ┼ёыш $x\in Inndiag(S)g$ ш $n$ эхўхЄэю, Єю $S_x\cong PSp_{n-1}(q)$;

    $(5)$ ┼ёыш $x\in Inndiag(S)g$, $n\geq 4$ ўхЄэю ш $p=2$, Єю яюфуЁєяяр $C_S(x)$ шчюьюЁЇэр ышсю $PSp_{n}(q)$, ышсю ЎхэЄЁрышчрЄюЁє эхъюЄюЁющ шэтюы■Ўшш т $PSp_{n}(q)$;

    $(6)$ ┼ёыш $x\in Inndiag(S)g$, $n\geq 4$ ўхЄэю ш $p>2$, Єю яюфуЁєяяр $S_x$ шчюьюЁЇэр юфэющ шч ёыхфє■∙шї уЁєяя: $PSp_{n}(q)$, $P\Omega^+_{n}(q)$, $P\Omega^-_{n}(q)$.
\end{lem}

\begin{lem}{\rm \cite[ЄхюЁхьр 4.5.1, яЁхфыюцхэш  2.5.12,\,4.9.1,\,4.9.2]{Maslova_GorLySo}}\label{AutPSU}
    ╧єёЄ№ $S = PSU_n(q)$, уфх $n\geq 3$, $q = p^m$ ш $(n,q)\not = (3,2)$, $x$ --- ¤ыхьхэЄ яЁюёЄюую яюЁ фър $r$ т $Aut(S)\setminus Inndiag(S)$ ш $S_x = O^{p'}(C_S(x))$. ╥юуфр ёяЁртхфышт√ ёыхфє■∙шх єЄтхЁцфхэш :

    $(1)$ $Aut(S)=Inndiag(S)\leftthreetimes \Phi$, уфх $Outdiag(S)\cong \mathbb Z_{(n,q+1)}$,
    $\Phi=\langle f \rangle\cong Aut(GF(q^2))\cong Z_{2m}$ --- уЁєяяр яюыхт√ї ртЄюьюЁЇшчьют уЁєяя√ $S$, $\Phi$ фхщёЄтєхЄ эр
    $Outdiag(S)$ Єръ цх, ъръ $Aut(GF(q^2))$ фхщёЄтєхЄ эр ьєы№ЄшяышърЄштэє■ яюфуЁєяяє шч $GF(q^2)$ Єюую цх яюЁ фър, ўЄю ш $Outdiag(S)$;

    $(2)$ ┼ёыш $r>2$, Єю $S_x\cong PSU_n(q^{\frac{1}{p}})$ ш $C_{Inndiag(S)}(x)\cong Inndiag(S_x)$;

    $(3)$ ┼ёыш $r=2$ ш $n$ эхўхЄэю, Єю $S_x\cong PSp_{n-1}(q)$;

    $(4)$ ┼ёыш $r=p=2$ ш $n$ ўхЄэю, Єю яюфуЁєяяр $C_S(x)$ шчюьюЁЇэр ышсю $PSp_{n}(q)$, ышсю ЎхэЄЁрышчрЄюЁє эхъюЄюЁющ шэтюы■Ўшш т $PSp_{n}(q)$;

    $(5)$ ┼ёыш $r=2$, $n$ ўхЄэю ш $p>2$, Єю яюфуЁєяяр $S_x$ шчюьюЁЇэр юфэющ шч ёыхфє■∙шї уЁєяя: $PSp_{n}(q)$, $P\Omega_{n}^+(q)$, $P\Omega_{n}^-(q)$.
\end{lem}

\section*{╩юэхўэ√х яюўЄш яЁюёЄ√х уЁєяя√, уЁрЇ√ ├Ё■эсхЁур--╩хухы  ъюЄюЁ√ї эх ёюфхЁцрЄ $3$-ъюъышъ}

┬ эрёЄю ∙хь Ёрчфхых фы  ърцфющ яюўЄш яЁюёЄющ уЁєяя√ $G$ ь√ яюърцхь, ўЄю ышсю $\Gamma(G)$ ёюфхЁцшЄ $3$-ъюъышъє, ышсю ёє∙хёЄтєхЄ ЁрчЁх°шьр  уЁєяяр $H$ Єрър , ўЄю $\Gamma(G)=\Gamma(H)$. ╥ръшь юсЁрчюь, ь√ фюърцхь ЄхюЁхьє \ref{AlmostSimple} ш шьяышърЎш■ $(1) \Rightarrow (3)$ т ЄхюЁхьх \ref{ASasS}. ╚ьяышърЎш  $(3) \Rightarrow (2)$ т ЄхюЁхьх \ref{ASasS} юўхтшфэр, р шьяышърЎш  $(2) \Rightarrow (1)$ ёыхфєхЄ шч \cite[ЄхюЁхьр~1]{Lucido2}.

╟рЇшъёшЁєхь юсючэрўхэшх, ъюЄюЁюх сєфхЄ фхщёЄтютрЄ№ фю ъюэЎр эрёЄю ∙хщ ЁрсюЄ√. ╧єёЄ№ $G$ --- яюўЄш яЁюёЄр  уЁєяяр ш $S = Soc(G)$.

\begin{lem}\label{Sporadic} ╧єёЄ№ $S$ --- ёяюЁрфшўхёър  уЁєяяр. ╥юуфр

$(1)$ $\Gamma(G)$ эх ёюфхЁцшЄ $3$-ъюъышъ, хёыш, ш Єюы№ъю хёыш $S \cong J_2${\rm;}

$(2)$ хёыш $S \cong J_2$, Єю ёє∙хёЄтєхЄ ЁрчЁх°шьр  уЁєяяр $H$ Єрър , ўЄю $\Gamma(G)=\Gamma(H)$.
\end{lem}

\proof ┬тшфє \cite{Atlas} шьххь $|G:S|\leq 2$.  ┼ёыш $S \not \cong J_2$, Єю ттшфє \cite[Єрсы. 1]{Maslova_VasVd2011} уЁрЇ $\Gamma(S)$ ёюфхЁцшЄ $3$-ъюъышъє, ёюёЄю ∙є■ шч эхўхЄэ√ї ўшёхы. ╧Ёшьхэхэшх ыхьь√ \ref{Maslova_nonadj} чртхЁ°рхЄ фюърчрЄхы№ёЄтю яєэъЄр $(1)$.

┼ёыш $S \cong J_2$, Єю ттшфє \cite{Atlas} шьххь $\pi(S)=\pi(Aut(S))=\{2,3,5\} \cup \{7\}$, ш ьэюцхёЄтр $\{2,3,5\}$ ш $\{7\}$  ты ■Єё  ъышърьш т $\Gamma(S)$. ╥ръшь юсЁрчюь, ттшфє ыхьь√ \ref{solvable} ёє∙хёЄтєхЄ ЁрчЁх°шьр  уЁєяяр $H$ Єрър , ўЄю $\Gamma(G)=\Gamma(H)$. ╦хььр фюърчрэр.

\begin{lem}\label{Alt} ╧єёЄ№ $S \cong A_n$, уфх $n\geq 5$. ╥юуфр

$(1)$ $\Gamma(G)$ эх ёюфхЁцшЄ $3$-ъюъышъ, хёыш, ш Єюы№ъю хёыш ышсю $A_6 < G \le Aut(A_6)$, ышсю $G \in \{S_5,S_8\}$, ышсю $S \in \{A_9, A_{10}, A_{12}\}${\rm;}

$(2)$ хёыш $\Gamma(G)$ эх ёюфхЁцшЄ $3$-ъюъышъ, Єю ёє∙хёЄтєхЄ ЁрчЁх°шьр  уЁєяяр $H$ Єрър , ўЄю $\Gamma(G)=\Gamma(H)$.
\end{lem}

\proof ┬тшфє \cite{Atlas} хёыш уЁєяяр $G$ шчюьюЁЇэр юфэющ шч уЁєяя: $A_5$, $A_6$, $A_7$, $S_7$, $A_8$, $A_{11}$, $S_{11}$, Єю $\Gamma(G)$ ёюфхЁцшЄ $3$-ъюышъє; хёыш уЁєяяр $G$ шчюьюЁЇэр  юфэющ шч уЁєяя: $S_5$, $S_6$, $M_{10}$, $PGL_2(9)$, $Aut(A_6)$, $S_8$, $A_9$, $S_9$, $A_{10}$, $S_{10}$, $A_{12}$, $S_{12}$, Єю $\Gamma(G)$ эх ёюфхЁцрЄ $3$-ъюъышъ, ш ыхуъю яЁютхЁшЄ№, ўЄю ттшфє ыхьь√ \ref{solvable} ёє∙хёЄтєхЄ ЁрчЁх°шьр  уЁєяяр $H$ Єрър , ўЄю $\Gamma(G)=\Gamma(H)$.

┼ёыш $13 \le n \le 17$, Єю ьэюцхёЄтю $\{7,11,13\}$ юсЁрчєхЄ $3$-ъюъышъє т $\Gamma(S)$ ттшфє \cite[яЁхфыюцхэшх~1.1]{Maslova_VasVd2005}. ┼ёыш $n\geq 18$, Єю яю \cite[ыхььр 1]{KM} т шэЄхЁтрых $((n + 1)/2,n)$ ыхцрЄ яю ъЁрщэхщ ьхЁх ЄЁш Ёрчышўэ√ї эхўхЄэ√ї яЁюёЄ√ї ўшёыр, ъюЄюЁ√х Єръцх юсЁрчє■Є $3$-ъюъышъє т $\Gamma(S)$ ттшфє \cite[яЁхфыюцхэшх~1.1]{Maslova_VasVd2005}. ╧Ёшьхэхэшх ыхьь√ \ref{Maslova_nonadj} чртхЁ°рхЄ фюърчрЄхы№ёЄтю ыхьь√.

\begin{lem}\label{PSL} ╧єёЄ№ $S \cong PSL_n(q)$, уфх $n\ge 2$, $q=p^m \ge 4$ ш $(n,q) \not = (2,2), (2,3)$. ╥юуфр

$(1)$ $\Gamma(G)$ эх ёюфхЁцшЄ $3$-ъюъышъ, хёыш, ш Єюы№ъю хёыш т√яюыэ хЄё  юфэю шч ёыхфє■∙шї єЄтхЁцфхэшщ{\rm:}

$\mbox{     }$$\mbox{     }$$\mbox{     }$$\mbox{     }$$\mbox{     }$$(i)$ уЁєяяр $G$ шчюьюЁЇэр юфэющ шч ёыхфє■∙шї уЁєяя{\rm:} $S_5$, $S_6$, $PGL_2(9)$, $M_{10}$, $Aut(A_6)$, $S_8$, $Aut(PSL_2(8))$,
$Aut(PSL_3(2))$, $PGL_3(4)\langle f \rangle$, $PGL_3(4)\langle g \rangle$, $Aut(PSL_3(4))$, $PSL_4(4)\langle f \rangle$, $PSL_4(4)\langle g \rangle$, $Aut(PSL_4(4))${\rm;}

$\mbox{     }$$\mbox{     }$$\mbox{     }$$\mbox{     }$$\mbox{     }$$(ii)$ $G$ --- уЁєяяр шч яєэъЄют $(3)-(9)$ ЄхюЁхь√ $\ref{AlmostSimple}${\rm;}

$(2)$ фы  ърцфющ уЁєяя√ $G$ шч яєэъЄр $(1)$ ыхьь√ уЁрЇ $\Gamma(G)$ Ёртхэ уЁрЇє ├Ё■эсхЁур--╩хухы  яюфїюф ∙хщ ЁрчЁх°шьющ уЁєяя√.

\end{lem}

\proof ╧єёЄ№ $n=2$. ╥юуфр ттшфє \cite[Єрсы.~5.1.A]{Maslova_KlLi} шьххь $|Out(S)|=(q-1,2)m$ ш $\pi(S)=\{p\} \cup R_1(q) \cup R_2(q)$.


╧єёЄ№ $q$ эхўхЄэю ш эх  ты хЄё  ўшёыюь тшфр $2^w \pm 1$. ╧єёЄ№ $r_1 \in R_{m}(p)$ ш $r_2 \in R_{2m}(p)$ эхўхЄэ√. ┬тшфє ьрыющ ЄхюЁхь√ ╘хЁьр ўшёыр $r_1$ ш $r_2$ эх фхы Є $|Out(S)|$. ┼ёыш $p$ эх фхышЄ шэфхъё $|G:S|$, Єю ттшфє ыхьь√ \ref{Maslova_nonadj} ўшёыр $\{p, r_1, r_2\}$ юсЁрчє■Є $3$-ъюъышъє т $\Gamma(G)$. ╧єёЄ№ $p$ фхышЄ шэфхъё $|G:S|$, Єюуфр $p$ фхышЄ $m$. ╧єёЄ№ $x$ --- ы■сющ ¤ыхьхэЄ яюЁ фър $p$, $y$ --- ы■сющ ¤ыхьхэЄ яюЁ фър $r_1$ ш $z$ --- ы■сющ ¤ыхьхэЄ яюЁ фър $r_2$ шч $G$. ╥юуфр $y$ ш $z$ ыхцрЄ т $Inndiag(S)$. ╟рьхЄшь, ўЄю т $Inndiag(S)$ эхЄ ¤ыхьхэЄют яюЁ фъют $pr_1$ ш $pr_2$. ┼ёыш $x \in G \setminus Inndiag(S)$, Єю ттшфє ыхьь√ \ref{AutPSL} шьххь $C_{Inndiag(S)}(x) \cong PGL_2(p^{m/p})$, ш ўшёыр $r_1$ ш $r_2$ эх фхы Є $|C_{Inndiag(S)}(x)|$. ╥ръшь юсЁрчюь, ўшёыр $\{p, r_1, r_2\}$ юсЁрчє■Є $3$-ъюъышъє т $\Gamma(G)$ ттшфє ыхьь√ \ref{Maslova_nonadj}.

╧єёЄ№ $q$ эхўхЄэю ш $q=2^w +\varepsilon1$, уфх $\varepsilon \in \{+,-\}$. ╥юуфр ттшфє ыхьь√ \ref{Maslova_Gerono} ышсю $q=9$ ш ёыєўрщ $S=PSL_2(9) \cong A_6$ ЁрёёьюЄЁхэ т ыхььх \ref{Alt}, ышсю $q=p$, $Aut(S) \cong PGL_2(p)$ ш $\pi(S)=\pi(Aut(S))$. ╟рьхЄшь, ўЄю $t(S)=3$ ттшфє \cite[Єрсы.~2]{Maslova_VasVd2011}. ┼ёыш $G$ --- уЁєяяр шч яєэъЄр $(3)$ ЄхюЁхь√ \ref{AlmostSimple}, Єю ттшфє \cite[яЁхфыюцхэш ~2.1,\,3.1,\,4.1]{Maslova_VasVd2005} ш \cite[ёыхфёЄтшх~2]{Maslova_But1} ьэюцхёЄтр $\{p\}$ ш $R_1(p)\cup R_2(p)$  ты ■Єё  ъышърьш т $\Gamma(G)$. ╥ръшь юсЁрчюь, ттшфє ыхьь√ \ref{solvable} ёє∙хёЄтєхЄ ЁрчЁх°шьр  уЁєяяр $H$ Єрър , ўЄю $\Gamma(G)=\Gamma(H)$.

╧єёЄ№ $q=2^m$. Cыєўрщ $S=PSL_2(4) \cong A_5$ ЁрёёьюЄЁхэ т ыхььх \ref{Alt}. ┬тшфє \cite{Atlas} уЁрЇ $\Gamma(PSL_2(8))$ ёюфхЁцшЄ $3$-ъюъышъє, р уЁрЇ $\Gamma(Aut(PSL_2(8)))$ эх ёюфхЁцшЄ $3$-ъюышъ. ╩Ёюьх Єюую, ттшфє \cite{Atlas} шьххь $\pi(Aut(PSL_2(8)))=\{2,3\} \cup \{7\}$, ш ўшёыр $2$ ш $3$ ёьхцэ√ т $\Gamma(Aut(PSL_2(8)))$. ╥ръшь юсЁрчюь, ттшфє ыхьь√ \ref{solvable} ёє∙хёЄтєхЄ ЁрчЁх°шьр  уЁєяяр $H$ Єрър , ўЄю $\Gamma(Aut(PSL_2(8)))=\Gamma(H)$.

╧єёЄ№ $q>8$, яюыюцшь $r_1=7$ ш $r_2=13$,  хёыш  $q=64$, ш $r_1 \in R_m(2)$ ш $r_2 \in R_{2m}(2)$ т яЁюЄштэюь ёыєўрх.
╟рьхЄшь, ўЄю ўшёыр $r_1$ ш $r_2$ эх фхы Є $m$ ттшфє ьрыющ ЄхюЁхь√ ╘хЁьр. ┼ёыш шэфхъё $|G:S|$ эхўхЄхэ, Єю ттшфє \cite[Єрсы.~2]{Maslova_VasVd2011}  ш ыхьь√ \ref{Maslova_nonadj} ўшёыр $\{2, r_1,r_2\}$ юсЁрчє■Є ъюъышъє т $\Gamma(G)$. ╧єёЄ№ $r \in \pi(|G:S|) \setminus \pi(S)$. ╥юуфр $r \not \in \{2,3\}$ ш $r$ фхышЄ $m$. ╧єёЄ№ $x$ --- ы■сющ ¤ыхьхэЄ яюЁ фър $r$, $y$ --- ы■сющ ¤ыхьхэЄ яюЁ фър $r_1$ ш $z$ --- ы■сющ ¤ыхьхэЄ яюЁ фър $r_2$ шч $G$. ╥юуфр $y$ ш $z$ ыхцрЄ т $Inndiag(S)$ ш $x \in G \setminus Inndiag(S)$. ┬тшфє ыхьь√ \ref{AutPSL} шьххь $C_{Inndiag(S)}(x) \cong PSL_2(p^{m/r})$ ш ўшёыр $r_1$ ш $r_2$ эх фхы Є $|C_{Inndiag(S)}(x)|$. ╥ръшь юсЁрчюь, ўшёыр $\{r, r_1, r_2\}$ юсЁрчє■Є $3$-ъюъышъє т $\Gamma(G)$. ╧єёЄ№ $\{2\} \subseteq \pi(|G:S|) \subset \pi(S)$. ┬тшфє \cite[яЁхфыюцхэш ~2.1,\,3.1,\,4.1]{Maslova_VasVd2005} ш ыхьь√ \ref{AutPSL} ьэюцхёЄтр $\{2\} \cup R_1(q)$ ш $R_2(q)$  ты ■Єё  ъышърьш т $\Gamma(G)$. ╥ръшь юсЁрчюь, хёыш $G$ --- уЁєяяр шч яєэъЄр $(4)$ ЄхюЁхь√ \ref{AlmostSimple}, Єю ттшфє ыхьь√ \ref{solvable} ёє∙хёЄтєхЄ ЁрчЁх°шьр  уЁєяяр $H$ Єрър , ўЄю $\Gamma(G)=\Gamma(H)$.

╧єёЄ№ $n=3$. ╥юуфр ттшфє \cite[Єрсы.~5.1.A]{Maslova_KlLi} шьххь $|Out(S)|=2(q-1,3)m$ ш $\pi(S)=\{p\}\cup R_1(q) \cup R_2(q) \cup R_3(q)$.

╧єёЄ№ $p$ эхўхЄэю ш $q+1$ эх  ты хЄё  ўшёыюь тшфр $2^w$. ┬тшфє \cite[Єрсы.~2]{Maslova_VasVd2011} т уЁрЇх $\Gamma(S)$ хёЄ№ ъюъышър $\{p, r_2, r_3\}$, уфх $2 \not = r_2 \in R_{2m}(p)$ ш $r_3 \in R_{3m}(p)$. ┬тшфє ьрыющ ЄхюЁхь√ ╘хЁьр ўшёыр $r_2$ ш $r_3$ эх фхы Є $|Out(S)|$.
┼ёыш $p$ эх фхышЄ шэфхъё $|G:S|$, Єю ўшёыр $\{p, r_2, r_3\}$ юсЁрчє■Є $3$-ъюъышъє т $\Gamma (G)$ ттшфє ыхьь√~\ref{Maslova_nonadj}. ╧єёЄ№ $p$ фхышЄ шэфхъё $|G:S|$. ╥юуфр $p$ фхышЄ $m$. ╧єёЄ№ $x$ --- ы■сющ ¤ыхьхэЄ яюЁ фър $p$, $y$ --- ы■сющ ¤ыхьхэЄ яюЁ фър $ r_2$ ш $z$ --- ы■сющ ¤ыхьхэЄ яюЁ фър $r_3$ шч $G$. ╥юуфр $y$ ш $z$ ыхцрЄ т $Inndiag(S)$. ╟рьхЄшь, ўЄю т $Inndiag(S)$ эхЄ ¤ыхьхэЄют яюЁ фъют $pr_2$ ш $pr_3$. ┼ёыш $x \in G \setminus Inndiag(S)$, Єю ттшфє ыхьь√ \ref{AutPSL} шьххь $C_{Inndiag(S)}(x)=PGL_3(p^{m/p})$, ш ўшёыр $r_2$ ш $r_3$ эх фхы Є $|C_{Inndiag(S)}(x)|$. ╥ръшь юсЁрчюь, ўшёыр $\{p, r_2, r_3\}$ юсЁрчє■Є $3$-ъюъышъє т $\Gamma(G)$.

╧Ёхфяюыюцшь, ўЄю $q=2^w - 1$. ╥юуфр ттшфє ыхьь√ \ref{Maslova_Gerono} шьххь $q=p$, $Aut(S) \cong PGL_3(q)\langle g \rangle$ ш $\pi(S)=\pi(Aut(S))$. ┼ёыш $(q-1)_3 \not =3$, Єю ттшфє \cite[яЁхфыюцхэш ~2.1,\,3.1,\,4.1]{Maslova_VasVd2005}  ьэюцхёЄтр $\{p\}\cup R_1(q) \cup R_2(q)$ ш $R_3(q)$  ты ■Єё  ъышърьш т $\Gamma(G)$. ╥ръшь юсЁрчюь, хёыш $G$ --- уЁєяяр шч яєэъЄр $(5)$ ЄхюЁхь√ \ref{AlmostSimple}, Єю ттшфє ыхьь√ \ref{solvable} ёє∙хёЄтєхЄ ЁрчЁх°шьр  уЁєяяр $H$ Єрър , ўЄю $\Gamma(G)=\Gamma(H)$. ╧єёЄ№ $(q-1)_3=3$. ╥юуфр ттшфє \cite[Єрсы.~2]{Maslova_VasVd2011} ьэюцхёЄтю $\{3,p,r_3\}$, уфх $r_3 \in R_3(p)$,  ты хЄё  $3$-ъюъышъющ т $\Gamma(S)$. ┼ёыш $\Gamma(G)$ эх ёюфхЁцшЄ $3$-ъюъышъ, Єю ттшфє ыхьь√ \ref{Maslova_nonadj} шьххь $Inndiag(S)\le G$. ╬сЁрЄэю, хёыш $Inndiag(S)\le G$, ттшфє \cite[яЁхфыюцхэш ~2.1,\,3.1,\,4.1]{Maslova_VasVd2005} ш \cite[ёыхфёЄтшх~2]{Maslova_But1} ьэюцхёЄтр $\{p\}\cup R_1(q) \cup R_2(q)$ ш $R_3(q)$  ты ■Єё  ъышърьш т $\Gamma(G)$. ╥ръшь юсЁрчюь, хёыш $G$ --- уЁєяяр шч яєэъЄр $(6)$ ЄхюЁхь√ \ref{AlmostSimple}, Єю ттшфє ыхьь√ \ref{solvable} ёє∙хёЄтєхЄ ЁрчЁх°шьр  уЁєяяр $H$ Єрър , ўЄю $\Gamma(G)=\Gamma(H)$.

╧єёЄ№ $p=2$. ╟рьхЄшь, ўЄю $PSL_3(2) \cong PSL_2(7)$, ш хёыш $q=4$, Єю ттшфє \cite[яЁхфыюцхэш ~2.1,\,3.1,\,4.1]{Maslova_VasVd2005} ш ыхьь√ \ref{AutPSL} уЁрЇ $\Gamma(G)$ эх ёюфхЁцшЄ $3$-ъюъышъ Єюуфр ш Єюы№ъю Єюуфр, ъюуфр $G \in \{PGL_3(4)\langle f \rangle, PGL_3(4)\langle g \rangle, Aut(PSL_3(4))\}$. ╩Ёюьх Єюую, $|\pi(Aut(PSL_3(4)))|=4$, ёыхфютрЄхы№эю, хёыш $\Gamma(G)$ эх ёюфхЁцшЄ $3$-ъюъышъ, Єю ттшфє ыхьь√ \ref{solvable} ёє∙хёЄтєхЄ ЁрчЁх°шьр  уЁєяяр $H$ Єрър , ўЄю $\Gamma(G)=\Gamma(H)$. ╧ю¤Єюьє ьюцэю ёўшЄрЄ№, ўЄю $m>2$.

 ┼ёыш $r \in \pi(|G:S|) \setminus \pi(S)$, Єю $r \not \in \{2,3\}$ ш $r$ фхышЄ $m$. ╧єёЄ№ $x$ --- ы■сющ ¤ыхьхэЄ яюЁ фър $r$, $y$ --- ы■сющ ¤ыхьхэЄ яюЁ фър $r_2 \in R_{2m}(2)$ ш $z$ --- ы■сющ ¤ыхьхэЄ яюЁ фър $r_3 \in R_{3m}(2)$ шч $G$. ╥юуфр $y$ ш $z$ ыхцрЄ т $Inndiag(S)$ ш $x \in G \setminus Inndiag(S)$. ┬тшфє ыхьь√ \ref{AutPSL} шьххь $C_{Inndiag(S)}(x)=PGL_3(2^{m/r})$ ш ўшёыр $r_2$ ш $r_3$ эх фхы Є $|C_{Inndiag(S)}(x)|$. ╥ръшь юсЁрчюь, ўшёыр $\{r, r_2, r_3\}$ юсЁрчє■Є $3$-ъюъышъє т $\Gamma(G)$.

╧єёЄ№ $\pi(G)=\pi(S)$ ш $\Gamma(G)$ эх ёюфхЁцшЄ $3$-ъюъышъ.

┼ёыш $(q-1)_3=3$, Єю $m \equiv 2,4 \pmod 6$ ш шч \cite[Єрсы.~2]{Maslova_VasVd2011} ёыхфєхЄ, ўЄю т уЁрЇх $\Gamma(S)$ хёЄ№ $3$-ъюъышър $\{3, r_2, r_3\}$, уфх $r_2 \in R_{2m}(2)$ ш $r_3 \in R_{3m}(2)$. ┬тшфє ьрыющ ЄхюЁхь√ ╘хЁьр ўшёыр $r_2$ ш $r_3$ эх фхы Є $|Out(S)|$. ╧юёъюы№ъє $3$ эх фхышЄ $m$, Єю ттшфє ыхьь√ \ref{Maslova_nonadj} шьххь $Inndiag(S)\le G$.

╧єёЄ№ $m$ эхўхЄэю. ╥юуфр $(q-1)_3=1$ ш $S=Inndiag(S)$. ┬тшфє \cite[Єрсы.~2]{Maslova_VasVd2011} т уЁрЇх $\Gamma(S)$ хёЄ№ $3$-ъюъышър $\{2, r_2, r_3\}$, уфх $r_3 \in R_{3m}(2)$, $r_2 \in R_{2m}(2)$, хёыш $m \not = 3$, ш $r_2 = 3$, хёыш $m=3$. ┬тшфє ыхьь√ \ref{Maslova_nonadj} шьххь $Inndiag(S)\langle g \rangle \le G$ .

╧єёЄ№ $m$ ўхЄэю ш $q \not = 16$.  ╨рёёьюЄЁшь $3$-ъюъышъє $\{2, r_2, r_3\}$ т $\Gamma(S)$, уфх $r_2 \in R_{2m}(2)$ ш $r_3 \in R_{3m}(2)$. $\Gamma(G)$ эх ёюфхЁцшЄ $3$-ъюъышъ, яю¤Єюьє ўшёыр $2$ ш $r_2$ ёьхцэ√ т $\Gamma(G)$ шыш ўшёыр $2$ ш $r_3$ ёьхцэ√ т $\Gamma(G)$. ┬ яхЁтюь ёыєўрх т $G \setminus Inndiag(S)$ эрщфхЄё  шэтюы■Ўш  $t$, яюЁ фюъ ЎхэЄЁрышчрЄюЁр ъюЄюЁющ т $Inndiag(S)$ фхышЄё  эр $r_2$. ┬тшфє ыхьь√ \ref{AutPSL} шьххь $t \in Inndiag(S)\langle g \rangle \setminus Inndiag(S)$, юЄъєфр $S\langle t\rangle \le G$. ╧Ёхфяюыюцшь, ўЄю ўшёыр $2$ ш $r_2$ эх ёьхцэ√ т $\Gamma(G)$, Єюуфр т $G \setminus Inndiag(S)$ эрщфхЄё  шэтюы■Ўш  $t_1$, яюЁ фюъ ЎхэЄЁрышчрЄюЁр ъюЄюЁющ т $Inndiag(S)$ фхышЄё  эр $r_3$. ┬тшфє ыхьь√ \ref{AutPSL} шьххь $t_1 \in Inndiag(S)f^{\frac{m}{2}}g$. ╨рёёьюЄЁшь $3$-ъюъышъє $\{2, r_2, r_3'\}$ т $\Gamma(S)$, уфх $r_3' \in R_{\frac{3m}{2}}(2)$. ╧юёъюы№ъє ўшёыр $2$ ш $r_2$ эх ёьхцэ√ т $\Gamma(G)$, Єюуфр т $G \setminus Inndiag(S)$ эрщфхЄё  шэтюы■Ўш  $t_2$, яюЁ фюъ ЎхэЄЁрышчрЄюЁр ъюЄюЁющ т $Inndiag(S)$ фхышЄё  эр $r_3'$. ┬тшфє ыхьь√ \ref{AutPSL} шьххь $t_2 \in Inndiag(S)\langle f^{\frac{m}{2}}\rangle \setminus Inndiag(S)$, юЄъєфр $t_1t_2 \in Inndiag(S)\langle g \rangle \setminus Inndiag(S)$ ш $S \langle t_1t_2 \rangle \le G$. ╩Ёюьх Єюую, ттшфє \cite[яЁхфыюцхэшх 4.9.2]{Maslova_GorLySo} ы■сющ ¤ыхьхэЄ яюЁ фър $2$ шч $Inndiag(S)\langle g \rangle \setminus Inndiag(S)$ ёюяЁ цхэ ё $g$ яюёЁхфёЄтюь ¤ыхьхэЄр шч $S$, юЄъєфр $S \langle g \rangle \le G$. ─ы  уЁєяя√ $S=PSL_3(16)$ рэрыюушўэ√щ Ёхчєы№ЄрЄ ьюцэю яюыєўшЄ№, яюыюцшт $r_2=17$, $r_3=13$ ш $r_3'=7$.

╧єёЄ№ $G$ --- уЁєяяр шч яєэъЄр $(7)$ шыш $(8)$ ЄхюЁхь√ \ref{AlmostSimple}. ╚ч \cite[яЁхфыюцхэш ~2.1,\,3.1,\,4.1]{Maslova_VasVd2005},  \cite[ёыхфёЄтшх~2]{Maslova_But1} ш \cite[яЁхфыюцхэшх 4.9.2]{Maslova_GorLySo} ёыхфєхЄ, ўЄю ьэюцхёЄтр $\{2\}\cup R_1(q) \cup R_2(q)$ ш $R_3(q)$  ты ■Єё  ъышърьш т $\Gamma(G)$. ╥ръшь юсЁрчюь, ттшфє ыхьь√ \ref{solvable} ёє∙хёЄтєхЄ ЁрчЁх°шьр  уЁєяяр $H$ Єрър , ўЄю $\Gamma(G)=\Gamma(H)$.

╧єёЄ№ $n=4$. ╥юуфр ттшфє \cite[Єрсы.~5.1.A]{Maslova_KlLi} шьххь $|Out(S)|=2(q-1,4)m$ ш $\pi(S)=\{p\}\cup R_1(q) \cup R_2(q) \cup R_3(q) \cup R_4(q)$.

╧єёЄ№ $p$ эхўхЄэю. ┬тшфє \cite[Єрсы.~2]{Maslova_VasVd2011} т уЁрЇх $\Gamma(S)$ хёЄ№ ъюъышър $\{p, r_3, r_4\}$, уфх $r_3 \in R_{3m}(p)$ ш $r_4 \in R_{4m}(p)$. ┬тшфє ьрыющ ЄхюЁхь√ ╘хЁьр ўшёыр $r_3$ ш $r_4$ эх фхы Є $|Out(S)|$. ┼ёыш $p$ эх фхышЄ шэфхъё $|G:S|$, Єю ўшёыр $\{p, r_3, r_4\}$ юсЁрчє■Є $3$-ъюъышъє т $\Gamma (G)$ ттшфє ыхьь√ \ref{Maslova_nonadj}. ╧єёЄ№ $p$ фхышЄ шэфхъё $|G:S|$. ╥юуфр $p$ фхышЄ $m$. ╧єёЄ№ $x$ --- ы■сющ ¤ыхьхэЄ яюЁ фър $p$, $y$ --- ы■сющ ¤ыхьхэЄ яюЁ фър $r_3$ ш $z$ --- ы■сющ ¤ыхьхэЄ яюЁ фър $r_4$ шч $G$. ╥юуфр $y$ ш $z$ ыхцрЄ т $Inndiag(S)$. ╟рьхЄшь, ўЄю т $Inndiag(S)$ эхЄ ¤ыхьхэЄют яюЁ фъют $pr_3$ ш $pr_4$. ┼ёыш $x \in G \setminus Inndiag(S)$, Єю ттшфє ыхьь√ \ref{AutPSL} шьххь $C_{Inndiag(S)}(x)=PGL_4(p^{m/p})$ ш ўшёыр $r_3$ ш $r_4$ эх фхы Є $|C_{Inndiag(S)}(x)|$. ╥ръшь юсЁрчюь, ўшёыр $\{p, r_3, r_4\}$ юсЁрчє■Є $3$-ъюъышъє т $\Gamma(G)$.

╧єёЄ№ $p=2$. ╥юуфр $S = Inndiag(S)$. ┼ёыш $q=2$, Єю, яюёъюы№ъє $PSL_4(2) \cong A_8$, ¤ЄюЄ ёыєўрщ ЁрёёьюЄЁхэ т ыхььх \ref{Alt}. ┼ёыш $q=4$, Єю ттшфє \cite[яЁхфыюцхэш ~2.1,\,3.1,\,4.1]{Maslova_VasVd2005} ш ыхьь√ \ref{AutPSL} уЁрЇ $\Gamma(G)$ эх ёюфхЁцшЄ $3$-ъюъышъ Єюуфр ш Єюы№ъю Єюуфр, ъюуфр $G \in \{PSL_4(4)\langle f \rangle, PGL_4(4)\langle g \rangle, Aut(PSL_4(4))\}$. ╩Ёюьх Єюую, $|\pi(Aut(PSL_4(4)))|=5$ ш $\Gamma(PSL_4(4))$ ёюфхЁцшЄ ЄЁхєуюы№эшъ, ёыхфютрЄхы№эю, $\Gamma(G)$ эх  ты хЄё  $5$-Ўшъыюь, юЄъєфр хёыш $\Gamma(G)$ эх ёюфхЁцшЄ $3$-ъюъышъ, Єю ттшфє ыхьь√ \ref{solvable} ёє∙хёЄтєхЄ ЁрчЁх°шьр  уЁєяяр $H$ Єрър , ўЄю $\Gamma(G)=\Gamma(H)$.  ╧ю¤Єюьє ьюцэю ёўшЄрЄ№, ўЄю $m>2$.

╧єёЄ№ $\Gamma(G)$ эх ёюфхЁцшЄ $3$-ъюъышъ. ╚ч \cite[Єрсы.~2]{Maslova_VasVd2011} ёыхфєхЄ, ўЄю т уЁрЇх $\Gamma(S)$ хёЄ№ $3$-ъюъышър $\{2, r_3, r_4\}$, уфх $r_3 \in R_{3m}(2)$ ш $r_m \in R_{4m}(2)$. ╧єёЄ№ $r \in \pi(|G:S|) \setminus \pi(S)$. ╥юуфр $r \not \in \{2,3\}$ ш $r$ фхышЄ $m$. ╧єёЄ№ $x$ --- ы■сющ ¤ыхьхэЄ яюЁ фър $r$, $y$ --- ы■сющ ¤ыхьхэЄ яюЁ фър $r_3$ ш $z$ --- ы■сющ ¤ыхьхэЄ яюЁ фър $r_4$ шч $G$. ╥юуфр $y$ ш $z$ ыхцрЄ т $Inndiag(S)$ ш $x \in G \setminus Inndiag(S)$. ┬тшфє ыхьь√ \ref{AutPSL} шьххь $C_{Inndiag(S)}(x)=PGL_4(2^{m/r})$ ш ўшёыр $r_3$ ш $r_4$ эх фхы Є $|C_{Inndiag(S)}(x)|$. ╥ръшь юсЁрчюь, ўшёыр $\{r, r_3, r_4\}$ юсЁрчє■Є $3$-ъюъышъє т $\Gamma(G)$, яЁюЄштюЁхўшх. ╧єёЄ№ $\pi(G)=\pi(S)$. ┼ёыш $m$ эхўхЄэю, Єю $Inndiag(S)\langle g \rangle \le G$ ттшфє ыхьь√ \ref{Maslova_nonadj}.
╧єёЄ№ $m$ ўхЄэю ш $q \not = 16$. ╥ръ ъръ $\Gamma(G)$ эх ёюфхЁцшЄ $3$-ъюъышъ, Єю ўшёыр $2$ ш $r_4$ ёьхцэ√ т $\Gamma(G)$ шыш ўшёыр $2$ ш $r_3$ ёьхцэ√ т $\Gamma(G)$. ┬ яхЁтюь ёыєўрх т $G \setminus Inndiag(S)$ эрщфхЄё  шэтюы■Ўш  $t$, яюЁ фюъ ЎхэЄЁрышчрЄюЁр ъюЄюЁющ т $Inndiag(S)$ фхышЄё  эр $r_4$. ┬тшфє ыхьь√ \ref{AutPSL} шьххь $t \in Inndiag(S)\langle g \rangle \setminus Inndiag(S)$, юЄъєфр $S\langle t\rangle \le G$. ╧Ёхфяюыюцшь, ўЄю ўшёыр $2$ ш $r_4$ эх ёьхцэ√ т $\Gamma(G)$, Єюуфр т $G \setminus Inndiag(S)$ эрщфхЄё  шэтюы■Ўш  $t_1$, яюЁ фюъ ЎхэЄЁрышчрЄюЁр ъюЄюЁющ т $Inndiag(S)$ фхышЄё  эр $r_3$. ┬тшфє ыхьь√ \ref{AutPSL} шьххь $t_1 \in Inndiag(S)f^{\frac{m}{2}}g$. ╚ч \cite[Єрсы.~2]{Maslova_VasVd2011} ёыхфєхЄ, ўЄю т уЁрЇх $\Gamma(S)$ хёЄ№ $3$-ъюъышър $\{2, r_3', r_4\}$, уфх $r_3' \in R_{\frac{3m}{2}}(2)$. ╧юёъюы№ъє ўшёыр $2$ ш $r_4$ эх ёьхцэ√ т $\Gamma(G)$, т $G \setminus Inndiag(S)$ эрщфхЄё  шэтюы■Ўш  $t_2$, яюЁ фюъ ЎхэЄЁрышчрЄюЁр ъюЄюЁющ т $Inndiag(S)$ фхышЄё  эр $r_3'$. ┬тшфє ыхьь√ \ref{AutPSL} шьххь $t_2 \in Inndiag(S)\langle f^{\frac{m}{2}} \rangle \setminus Inndiag(S)$, юЄъєфр $t_1t_2 \in Inndiag(S)\langle g \rangle \setminus Inndiag(S)$ ш $S \langle t_1t_2 \rangle \le G$, яю¤Єюьє $Inndiag(S)\langle g \rangle \le G$. ─ы  уЁєяя√ $S=PSL_4(16)$ рэрыюушўэ√щ Ёхчєы№ЄрЄ ьюцэю яюыєўшЄ№, яюыюцшт $r_4=257$, $r_3=13$ ш $r_3'=7$.

╧єёЄ№ $G$ --- уЁєяяр шч яєэъЄр $(9)$ ЄхюЁхь√ \ref{AlmostSimple}. ╚ч \cite[яЁхфыюцхэш ~2.1,\,3.1,\,4.1]{Maslova_VasVd2005} ш \cite[яЁхфыюцхэшх 4.9.2]{Maslova_GorLySo} ёыхфєхЄ, ўЄю ьэюцхёЄтр $\{2\}\cup R_2(q) \cup R_4(q)$ ш $R_1(q) \cup R_3(q)$  ты ■Єё  ъышърьш т $\Gamma(G)$. ╥ръшь юсЁрчюь, ттшфє ыхьь√ \ref{solvable} ёє∙хёЄтєхЄ ЁрчЁх°шьр  уЁєяяр $H$ Єрър , ўЄю $\Gamma(G)=\Gamma(H)$.

╧єёЄ№ $n \ge 5$. ╥юуфр ттшфє \cite[Єрсы.~5.1.A]{Maslova_KlLi} шьххь $|Out(S)|=2(q-1,n)m$. ╨рёёьюЄЁшь ьэюцхёЄтю $\{a,b,c\}$, уфх $a \in R_{mn}(p) \subseteq R_{n}(q)$, $b \in R_{m(n-1)}(p) \subseteq R_{n-1}(q)$ ш $c \in R_{m(n-2)}(p) \subseteq R_{n-2}(q)$, хёыш $(m,n,p) \not = (2,5,2),(1,6,2),(1,7,2),(1,8,2)$; $\{a,b,c\}=\{7,17,31\}$, хёыш $(m,n,p) = (2,5,2)$; $\{a,b,c\}=\{5,7,31\}$, хёыш $(m,n,p) = (1,6,2)$; $\{a,b,c\}=\{5,31, 127\}$, хёыш $(m,n,p) = (1,7,2),(1,8,2)$. ┬тшфє ьрыющ ЄхюЁхь√ ╘хЁьр ўшёыр $a$, $b$ ш $c$ эх фхы Є $|Out(S)|$. ╚ч \cite[Єрсы.~2]{Maslova_VasVd2011} ёыхфєхЄ, ўЄю ўшёыр $\{a,b,c\}$ юсЁрчє■Є ъюъышъє т $\Gamma(S)$, яю¤Єюьє юэш юсЁрчє■Є $3$-ъюъышъє т $\Gamma(G)$ ттшфє ыхьь√ \ref{Maslova_nonadj}. ╦хььр фюърчрэр.

\begin{lem}\label{PSU} ╧єёЄ№ $S \cong PSU_n(q)$, уфх $n\ge 3$, $q=p^m \ge 2$ ш $(n,q) \not = (3,2)$. ╥юуфр

$(1)$ $\Gamma(G)$ эх ёюфхЁцшЄ $3$-ъюъышъ, хёыш, ш Єюы№ъю хёыш т√яюыэ хЄё  юфэю шч ёыхфє■∙шї єЄтхЁцфхэшщ{\rm:}

$\mbox{     }$$\mbox{     }$$\mbox{     }$$\mbox{     }$$\mbox{     }$$(i)$ $S$ шчюьюЁЇэр юфэющ шч уЁєяя $PSU_3(9)$ шыш $PSU_4(2)${\rm;}

$\mbox{     }$$\mbox{     }$$\mbox{     }$$\mbox{     }$$\mbox{     }$$(ii)$ $G \cong Aut(PSU_5(2))${\rm;}

$\mbox{     }$$\mbox{     }$$\mbox{     }$$\mbox{     }$$\mbox{     }$$(iii)$ $G$ --- уЁєяяр шч яєэъЄют $(10)-(14)$ ЄхюЁхь√ $\ref{AlmostSimple}${\rm;}

$(2)$ фы  ърцфющ уЁєяя√ $G$ шч яєэъЄр $(1)$ ыхьь√ уЁрЇ $\Gamma(G)$ Ёртхэ уЁрЇє ├Ё■эсхЁур--╩хухы  яюфїюф ∙хщ ЁрчЁх°шьющ уЁєяя√.

\end{lem}

\proof ╧єёЄ№ $n=3$. ╥юуфр ттшфє \cite[Єрсы.~5.1.A]{Maslova_KlLi} шьххь $|Out(S)|=2m(q+1,3)$ ш $\pi(S)=\{p\}\cup R_1(q) \cup R_2(q) \cup R_6(q)$.

╧єёЄ№ $p$ эхўхЄэю ш $q-1$ эх  ты хЄё  ўшёыюь тшфр $2^w$. ┬тшфє \cite[Єрсы.~2]{Maslova_VasVd2011} т уЁрЇх $\Gamma(S)$ хёЄ№ ъюъышър $\{p, 3, r_1, r_6\}$, уфх $2 \not= r_1 \in R_m(p)$ ш $r_6 \in R_{6m}(p)$. ┬тшфє ьрыющ ЄхюЁхь√ ╘хЁьр ўшёыр $r_1$ ш $r_6$ эх фхы Є $|Out(S)|$.
┼ёыш $p$ эх фхышЄ шэфхъё $|G:S|$, Єю ўшёыр $\{p, r_1, r_6\}$ юсЁрчє■Є $3$-ъюъышъє т $\Gamma (G)$ ттшфє ыхьь√ \ref{Maslova_nonadj}. ╧єёЄ№ $p$ фхышЄ шэфхъё $|G:S|$. ╥юуфр $p$ фхышЄ $m$. ╧єёЄ№ $x$ --- ы■сющ ¤ыхьхэЄ яюЁ фър $p$, $y$ --- ы■сющ ¤ыхьхэЄ яюЁ фър $r_1$ ш $z$ --- ы■сющ ¤ыхьхэЄ яюЁ фър $r_6$ шч $G$. ╥юуфр $y$ ш $z$ ыхцрЄ т $Inndiag(S)$. ╟рьхЄшь, ўЄю т $Inndiag(S)$ эхЄ ¤ыхьхэЄют яюЁ фъют $pr_1$ ш $pr_6$. ┼ёыш $x \in G \setminus Inndiag(S)$, Єю ттшфє ыхьь√ \ref{AutPSL} шьххь $C_{Inndiag(S)}(x)=PGU_3(p^{m/p})$, ш ўшёыр $r_1$ ш $r_6$ эх фхы Є $|C_{Inndiag(S)}(x)|$. ╥ръшь юсЁрчюь, ўшёыр $\{p, r_1, r_6\}$ юсЁрчє■Є $3$-ъюъышъє т $\Gamma(G)$.

╧Ёхфяюыюцшь, ўЄю $q=2^w + 1$. ╥юуфр ттшфє ыхьь√ \ref{Maslova_Gerono} шьххь ышсю $q=9$, ышсю $q=p$. ┬ юсюшї ёыєўр ї $\pi(S)=\pi(Aut(S))$. ┼ёыш $(q+1)_3 \not =3$, Єю ттшфє \cite[яЁхфыюцхэш ~2.2,\,3.1,\,4.2]{Maslova_VasVd2005} ьэюцхёЄтр $\{p\}\cup R_1(q) \cup R_2(q)$ ш $R_6(q)$  ты ■Єё  ъышърьш т $\Gamma(G)$. ╥ръшь юсЁрчюь, хёыш $S \cong PSU_3(9)$ шыш $G$ --- уЁєяяр шч яєэъЄр $(10)$ ЄхюЁхь√ \ref{AlmostSimple}, Єю ттшфє ыхьь√ \ref{solvable} ёє∙хёЄтєхЄ ЁрчЁх°шьр  уЁєяяр $H$ Єрър , ўЄю $\Gamma(G)=\Gamma(H)$. ╧єёЄ№ $(q+1)_3=3$. ╥юуфр ттшфє \cite[Єрсы.~2]{Maslova_VasVd2011} ьэюцхёЄтю  $\{3,p,r_6\}$, уфх $r_6 \in R_6(p)$,  ты хЄё  ъюъышъющ т $\Gamma(S)$. ┼ёыш $\Gamma(G)$ эх ёюфхЁцшЄ $3$-ъюъышъ, Єю ттшфє ыхьь√ \ref{Maslova_nonadj} шьххь $3$ фхышЄ $|G/S|$, юЄъєфр ёыхфєхЄ, ўЄю $Inndiag(S)\le G$. ╬сЁрЄэю, хёыш $Inndiag(S)\le G$, ттшфє \cite[яЁхфыюцхэш ~2.1,\,3.1,\,4.1]{Maslova_VasVd2005} ш \cite[ёыхфёЄтшх~2]{Maslova_But1} ьэюцхёЄтр $\{p\}\cup R_1(q) \cup R_2(q)$ ш $R_6(q)$  ты ■Єё  ъышърьш т $\Gamma(G)$. ╥ръшь юсЁрчюь, хёыш $G$ --- уЁєяяр шч яєэъЄр $(11)$ ЄхюЁхь√ \ref{AlmostSimple}, Єю ттшфє ыхьь√ \ref{solvable} ёє∙хёЄтєхЄ ЁрчЁх°шьр  уЁєяяр $H$ Єрър , ўЄю $\Gamma(G)=\Gamma(H)$.

╧єёЄ№ $p=2$. ┼ёыш $(q+1)_3=3$, Єю $m \equiv \pm 1 \pmod 6$. ╚ч \cite[Єрсы.~2]{Maslova_VasVd2011} ёыхфєхЄ, ўЄю т уЁрЇх $\Gamma(S)$ хёЄ№
$3$-ъюъышър $\{3, r_1, r_6\}$, уфх $r_1 \in R_{m}(2)$ ш $r_6 \in R_{6m}(2)$. ┬тшфє ьрыющ ЄхюЁхь√ ╘хЁьр ўшёыр $r_1$ ш $r_6$ эх фхы Є $|Out(S)|$. ╧юёъюы№ъє $3$ эх фхышЄ $m$, Єю хёыш $\Gamma(G)$ эх ёюфхЁцшЄ $3$-ъюъышъ, ттшфє ыхьь√ \ref{Maslova_nonadj} шьххь $Inndiag(S)\le G$.

╧єёЄ№ $\Gamma(G)$ эх ёюфхЁцшЄ $3$-ъюъышъ. ╧юыюцшь $r_6 \in R_{6m}(2)$ ш $r_1 \in R_{m}(2)$. ┬тшфє ьрыющ ЄхюЁхь√ ╘хЁьр ўшёыр $r_6$ ш $r_1$ эх фхы Є $|Out(S)|$. ╧Ёхфяюыюцшь, ўЄю $r \in \pi(|G:S|) \setminus \pi(S)$, Єюуфр $r \not \in \{2,3\}$ ш $r$ фхышЄ $m$. ╧єёЄ№ $x$ --- ы■сющ ¤ыхьхэЄ яюЁ фър $r$, $y$ --- ы■сющ ¤ыхьхэЄ яюЁ фър $r_1$ ш $z$ --- ы■сющ ¤ыхьхэЄ яюЁ фър $r_6$ шч $G$. ╥юуфр $y$ ш $z$ ыхцрЄ т $Inndiag(S)$ ш $x \in G \setminus Inndiag(S)$. ┬тшфє ыхьь√ \ref{AutPSL} шьххь $C_{Inndiag(S)}(x)=PGU_3(2^{m/r})$ ш ўшёыр $r_1$ ш $r_6$ эх фхы Є $|C_{Inndiag(S)}(x)|$. ╥ръшь юсЁрчюь, ўшёыр $\{r, r_1, r_6\}$ юсЁрчє■Є $3$-ъюъышъє т $\Gamma(G)$, яЁюЄштюЁхўшх. ╧єёЄ№ $\pi(G)=\pi(S)$ ш шэфхъё $|G:S|$ эхўхЄхэ. ╥юуфр ўшёыр $\{2, r_1, r_6\}$ юсЁрчє■Є $3$-ъюъышъє т $\Gamma(G)$ ттшфє \cite[Єрсы.~2]{Maslova_VasVd2011} ш ыхьь√ \ref{Maslova_nonadj}. ╧єёЄ№ $G$ --- уЁєяяр шч яєэъЄр $(12)$ шыш $(13)$ ЄхюЁхь√ \ref{AlmostSimple}. ┬тшфє \cite[яЁхфыюцхэш ~2.2,\,3.1,\,4.2]{Maslova_VasVd2005} ш \cite[ёыхфёЄтшх~2]{Maslova_But1} ьэюцхёЄтр $\{2\}\cup R_1(q) \cup R_2(q)$ ш $R_6(q)$  ты ■Єё  ъышърьш т $\Gamma(G)$. ╥ръшь юсЁрчюь, ттшфє ыхьь√ \ref{solvable} ёє∙хёЄтєхЄ ЁрчЁх°шьр  уЁєяяр $H$ Єрър , ўЄю $\Gamma(G)=\Gamma(H)$.

╧єёЄ№ $n=4$. ╥юуфр ттшфє \cite[Єрсы.~5.1.A]{Maslova_KlLi} шьххь $|Out(S)|=2(q+1,4)m$ ш $\pi(S)=\{p\}\cup R_1(q) \cup R_2(q) \cup R_4(q) \cup R_6(q)$.

╧єёЄ№ $p$ эхўхЄэю. ┬тшфє \cite[Єрсы.~2]{Maslova_VasVd2011} т уЁрЇх $\Gamma(S)$ хёЄ№ ъюъышър $\{p, r_4, r_6\}$, уфх $r_4 \in R_{4m}(p)$ ш $r_6 \in R_{6m}(p)$. ┬тшфє ьрыющ ЄхюЁхь√ ╘хЁьр ўшёыр $r_4$ ш $r_6$ эх фхы Є $|Out(S)|$. ┼ёыш $p$ эх фхышЄ шэфхъё $|G:S|$, Єю ўшёыр $\{p, r_4, r_6\}$ юсЁрчє■Є $3$-ъюъышъє т $\Gamma (G)$ ттшфє ыхьь√ \ref{Maslova_nonadj}. ╧єёЄ№ $p$ фхышЄ шэфхъё $|G:S|$, Єюуфр $p$ фхышЄ $m$. ╧єёЄ№ $x$ --- ы■сющ ¤ыхьхэЄ яюЁ фър $p$, $y$ --- ы■сющ ¤ыхьхэЄ яюЁ фър $r_4$ ш $z$ --- ы■сющ ¤ыхьхэЄ яюЁ фър $r_6$ шч $G$. ╥юуфр $y$ ш $z$ ыхцрЄ т $Inndiag(S)$. ╟рьхЄшь, ўЄю т $Inndiag(S)$ эхЄ ¤ыхьхэЄют яюЁ фъют $pr_4$ ш $pr_6$. ┼ёыш $x \in G \setminus Inndiag(S)$, Єю ттшфє ыхьь√ \ref{AutPSL} шьххь $C_{Inndiag(S)}(x)=PGU_4(p^{m/p})$ ш ўшёыр $r_4$ ш $r_6$ эх фхы Є $|C_{Inndiag(S)}(x)|$. ╥ръшь юсЁрчюь, ўшёыр $\{p, r_4, r_6\}$ юсЁрчє■Є $3$-ъюъышъє т $\Gamma(G)$.

╧єёЄ№ $p=2$. ╥юуфр $S = Inndiag(S)$. ┬тшфє \cite{Atlas} шьххь $\pi(PSU_4(2))=\pi(Aut(PSU_4(2)))$ ш $\pi(PSU_4(2))=\{2,3\} \cup \{5\}$, яЁшўхь ўшёыр $2$ ш $3$ ёьхцэ√ т $\Gamma(PSU_4(2))$. ╥ръшь юсЁрчюь, ттшфє ыхьь√ \ref{solvable} хёыш $S=PSU_4(2)$, Єю ёє∙хёЄтєхЄ ЁрчЁх°шьр  уЁєяяр $H$ Єрър , ўЄю $\Gamma(S)=\Gamma(H)$. ╧ю¤Єюьє ьюцхь ёўшЄрЄ№, ўЄю $q \ge 4$.
╚ч \cite[Єрсы.~2]{Maslova_VasVd2011} ёыхфєхЄ, ўЄю т уЁрЇх $\Gamma(S)$ хёЄ№ $3$-ъюъышър $\{2, r_4, r_6\}$, уфх $r_4 \in R_{4m}(2)$ ш $r_6 \in R_{6m}(2)$. ┬тшфє ьрыющ ЄхюЁхь√ ╘хЁьр ўшёыр $r_4$ ш $r_6$ эх фхы Є $|Out(S)|$.
╧єёЄ№ $r \in \pi(|G:S|) \setminus \pi(S)$. ╥юуфр $r \not \in \{2,3\}$ ш $r$ фхышЄ $m$. ╧єёЄ№ $x$ --- ы■сющ ¤ыхьхэЄ яюЁ фър $r$, $y$ --- ы■сющ ¤ыхьхэЄ яюЁ фър $r_4$ ш $z$ --- ы■сющ ¤ыхьхэЄ яюЁ фър $r_6$ шч $G$. ╥юуфр $y$ ш $z$ ыхцрЄ т $Inndiag(S)$ ш $x \in G \setminus Inndiag(S)$. ┬тшфє ыхьь√ \ref{AutPSU} шьххь $C_{Inndiag(S)}(x)=PGU_4(2^{m/r})$ ш ўшёыр $r_4$ ш $r_6$ эх фхы Є $|C_{Inndiag(S)}(x)|$. ╥ръшь юсЁрчюь, ўшёыр $\{r, r_4, r_6\}$ юсЁрчє■Є $3$-ъюъышъє т $\Gamma(G)$. ╧єёЄ№ $\pi(G)=\pi(S)$ ш шэфхъё $|G:S|$ эхўхЄхэ. ╥юуфр ўшёыр $\{2, r_4, r_6\}$ юсЁрчє■Є $3$-ъюъышъє т $\Gamma(G)$ ттшфє \cite[Єрсы.~2]{Maslova_VasVd2011} ш ыхьь√ \ref{Maslova_nonadj}. ╧єёЄ№ $\pi(G)=\pi(S)$ ш шэфхъё $|G:S|$ ўхЄхэ, Єюуфр $m$ ўхЄэю ш $G$ ёюфхЁцшЄ шэтюы■Ўш■ $t \in G \setminus Inndiag(S)$, ёыхфютрЄхы№эю, $S \langle t \rangle \le G$, юЄъєфр т $G$ хёЄ№ шэтюы■Ўш  $f^{\frac{m}{2}}$. ┼ёыш $G$ --- уЁєяяр шч яєэъЄр $(14)$ ЄхюЁхь√ \ref{AlmostSimple}, Єю ттшфє \cite[яЁхфыюцхэш ~2.2,\,3.1,\,4.2]{Maslova_VasVd2005} ш \cite[2.5.13, 4.9.2]{Maslova_GorLySo} ьэюцхёЄтр $\{2\}\cup R_1(q) \cup R_4(q)$ ш $R_2(q) \cup R_6(q)$  ты ■Єё  ъышърьш т $\Gamma(G)$, яю¤Єюьє ттшфє ыхьь√ \ref{solvable} ёє∙хёЄтєхЄ ЁрчЁх°шьр  уЁєяяр $H$ Єрър , ўЄю $\Gamma(G)=\Gamma(H)$.

╧єёЄ№ $n = 5$. ╥юуфр ттшфє \cite[Єрсы.~5.1.A]{Maslova_KlLi} шьххь $|Out(S)|=2(q+1,5)m$. ┼ёыш $q=2$, Єю ттшфє \cite[яЁхфыюцхэш ~2.2,\,3.1,\,4.2]{Maslova_VasVd2005} ш \cite[2.5.13, 4.9.2]{Maslova_GorLySo} уЁрЇ $\Gamma(PSU_5(2))$ ёюфхЁцшЄ $3$-ъюъышъє $\{2,5,11\}$, ш $\pi(Aut(PSU_5(2)))=\{2,3,5\} \cup \{11\}$, уфх $\{2,3,5\}$ --- ъышър т $\Gamma(Aut(PSU_5(2)))$. ╤ыхфютрЄхы№эю, хёыш $G \cong Aut(PSU_5(2))$, Єю ттшфє ыхьь√ \ref{solvable} ёє∙хёЄтєхЄ ЁрчЁх°шьр  уЁєяяр $H$ Єрър , ўЄю $\Gamma(G)=\Gamma(H)$. ╧ю¤Єюьє ьюцэю ёўшЄрЄ№, ўЄю $m>1$. ╨рёёьюЄЁшь ьэюцхёЄтю $\{r_4,r_6,r_{10}\}$, уфх $r_4 \in R_{4m}(p) \subseteq R_{4}(q)$, $r_6 \in R_{6m}(p) \subseteq R_{6}(q)$ ш $r_{10} \in R_{10m}(p) \subseteq R_{10}(q)$. ┬тшфє ьрыющ ЄхюЁхь√ ╘хЁьр ўшёыр $r_4$, $r_6$ ш $r_{10}$ эх фхы Є $|Out(S)|$. ╚ч \cite[Єрсы.~2]{Maslova_VasVd2011} ёыхфєхЄ, ўЄю ўшёыр $\{r_4,r_6,r_{10}\}$ юсЁрчє■Є $3$-ъюъышъє т $\Gamma(S)$, яю¤Єюьє юэш юсЁрчє■Є $3$-ъюъышъє т $\Gamma(G)$ ттшфє ыхьь√ \ref{Maslova_nonadj}.

╧єёЄ№ $n = 6$. ╥юуфр ттшфє \cite[Єрсы.~5.1.A]{Maslova_KlLi} шьххь $|Out(S)|=2(q+1,6)m$. ╨рёёьюЄЁшь ьэюцхёЄтю $\{r_3,r_4,r_{10}\}$, уфх $r_3 \in R_{3m}(p)$, $r_4 \in R_{4m}(p)$ ш $r_{10} \in R_{10m}(p)$, хёыш $(m,p) \not = (2,2)$, ш  $r_3 =13$, $r_4 =17$ ш $r_{10} =41$, хёыш $(m,p) = (2,2)$. ┬тшфє ьрыющ ЄхюЁхь√ ╘хЁьр ўшёыр $r_3$, $r_4$ ш $r_{10}$ эх фхы Є $|Out(S)|$. ╚ч \cite[Єрсы.~2]{Maslova_VasVd2011} ёыхфєхЄ, ўЄю ўшёыр $\{r_3,r_4,r_{10}\}$ юсЁрчє■Є $3$-ъюъышъє т $\Gamma(S)$, яю¤Єюьє юэш юсЁрчє■Є $3$-ъюъышъє т $\Gamma(G)$ ттшфє ыхьь√ \ref{Maslova_nonadj}.

╧єёЄ№ $n \ge 7$. ╥юуфр ттшфє \cite[Єрсы.~5.1.A]{Maslova_KlLi} шьххь $|Out(S)|=2(q+1,n)m$. ╨рёёьюЄЁшь ьэюцхёЄтю $\{a,b,c\}$, уфх $a \in R_{2mn}(p)$, $b \in R_{2m(n-2)}(p)$ ш $c \in R_{m(n-2+\varepsilon1)}(p)$, хёыш $\varepsilon \in \{+,-\}$ ш $n \equiv \varepsilon 1 \pmod 4$; $a \in R_{2m(n-1)}(p)$, $b \in R_{2m(n-3)}(p)$ ш $c \in R_{mn}(p)$, хёыш $n \equiv 0 \pmod 4$; $a \in R_{2m(n-1)}(p)$, $b \in R_{2m(n-3)}(p)$ ш $c \in R_{m(n-2)}(p)$, хёыш $n \equiv 2 \pmod 4$. ┬тшфє ьрыющ ЄхюЁхь√ ╘хЁьр ўшёыр $a$, $b$ ш $c$ эх фхы Є $|Out(S)|$. ╚ч \cite[Єрсы.~2]{Maslova_VasVd2011} ёыхфєхЄ, ўЄю ўшёыр $\{a,b,c\}$ юсЁрчє■Є $3$-ъюъышъє т $\Gamma(S)$, яю¤Єюьє юэш юсЁрчє■Є $3$-ъюъышъє т $\Gamma(G)$ ттшфє ыхьь√ \ref{Maslova_nonadj}. ╦хььр фюърчрэр.

\begin{lem}\label{PSp} ╧єёЄ№ $S \cong PSp_n(q)$, уфх $n \ge 4$ ўхЄэю, $q=p^m \ge 2$ ш $(n,q) \not = (4,2)$. ╥юуфр

$(1)$ $\Gamma(G)$ эх ёюфхЁцшЄ $3$-ъюъышъ, хёыш, ш Єюы№ъю хёыш т√яюыэ хЄё  юфэю шч ёыхфє■∙шї єЄтхЁцфхэшщ{\rm:}

$\mbox{     }$$\mbox{     }$$\mbox{     }$$\mbox{     }$$\mbox{     }$$(i)$ $S \cong PSp_6(2)${\rm;}

$\mbox{     }$$\mbox{     }$$\mbox{     }$$\mbox{     }$$\mbox{     }$$(ii)$ $G$ --- уЁєяяр шч яєэъЄр $(15)$ ЄхюЁхь√ $\ref{AlmostSimple}${\rm;}

$(2)$ фы  ърцфющ уЁєяя√ $G$ шч яєэъЄр $(1)$ ыхьь√ уЁрЇ $\Gamma(G)$ Ёртхэ уЁрЇє ├Ё■эсхЁур--╩хухы  яюфїюф ∙хщ ЁрчЁх°шьющ уЁєяя√.

\end{lem}

\proof ╨рёёьюЄЁшь ёыєўрщ $n=4$. ╧єёЄ№ $r \in \pi(|G:S|) \setminus \pi(S)$, Єюуфр $r \not \in \{2,3\}$ ш $r$ фхышЄ $m$. ╧єёЄ№ $x$ --- ы■сющ ¤ыхьхэЄ яюЁ фър $r$, $y$ --- ы■сющ ¤ыхьхэЄ яюЁ фър $r_2 \in R_{2m}(p)$ ш $z$ --- ы■сющ ¤ыхьхэЄ яюЁ фър $r_4 \in R_{4m}(p)$ шч $G$. ╥юуфр $y$ ш $z$ ыхцрЄ т $Inndiag(S)$ ш $x \in G \setminus Inndiag(S)$. ┬тшфє \cite[яЁхфыюцхэш  4.9.1]{Maslova_GorLySo} шьххь $C_{Inndiag(S)}(x)=PSp_4(p^{m/r})$ ш ўшёыр $r_2$ ш $r_4$ эх фхы Є $|C_{Inndiag(S)}(x)|$. ╥ръшь юсЁрчюь, ўшёыр $\{r, r_2, r_4\}$ юсЁрчє■Є $3$-ъюъышъє т $\Gamma(G)$. ╧єёЄ№ $\pi(G)=\pi(S)$.  ┬тшфє \cite[ЄхюЁхьр~1]{ZinMaz} $\Gamma(S)$ ёюёЄюшЄ шч фтєї эхяхЁхёхър■∙шїё  ъышъ, ёыхфютрЄхы№эю ттшфє ыхьь√ \ref{solvable} ёє∙хёЄтєхЄ ЁрчЁх°шьр  уЁєяяр $H$ Єрър , ўЄю $\Gamma(G)=\Gamma(H)$.

╧єёЄ№ $n \ge 6$. ╥юуфр ттшфє \cite[Єрсы.~5.1.A]{Maslova_KlLi} шьххь $|Out(S)|=2(q-1,2)m$. ┼ёыш $(n,q)=(6,2)$, Єю $\pi(S)=\pi(Aut(S))$ ш  $\Gamma(S)$ ёюёЄюшЄ шч фтєї эхяхЁхёхър■∙шїё  ъышъ ттшфє \cite[ЄхюЁхьр~1]{ZinMaz}, ёыхфютрЄхы№эю ттшфє ыхьь√ \ref{solvable} ёє∙хёЄтєхЄ ЁрчЁх°шьр  уЁєяяр $H$ Єрър , ўЄю $\Gamma(G)=\Gamma(H)$. ┼ёыш $n=6$ ш $q>2$, Єю ЁрёёьюЄЁшь ьэюцхёЄтю $\{r_3,r_4,r_6\}$, уфх $r_3 \in R_{3m}(p)$, $r_4 \in R_{4m}(p)$ ш $r_6 \in R_{6m}(p)$, хёыш $(m,p) \not = (2,2)$, ш $r_3=7$, $r_4 =17$ ш $r_6=13$, хёыш $(m,p)  = (2,2)$. ╟рьхЄшь, ўЄю ўшёыр $r_3$, $r_4$ ш $r_6$ эхўхЄэ√ ш юЄышўэ√ юЄ $3$. ┬тшфє ьрыющ ЄхюЁхь√ ╘хЁьр ўшёыр $r_3$, $r_4$ ш $r_6$ эх фхы Є $|Out(S)|$. ╚ч \cite[Єрсы.~3]{Maslova_VasVd2011} ёыхфєхЄ, ўЄю ўшёыр $\{r_3,r_4,r_6\}$ юсЁрчє■Є $3$-ъюъышъє т $\Gamma(S)$, яю¤Єюьє юэш юсЁрчє■Є $3$-ъюъышъє т $\Gamma(G)$ ттшфє ыхьь√ \ref{Maslova_nonadj}.

╧єёЄ№ $n \ge 8$. ╨рёёьюЄЁшь ьэюцхёЄтю $\{a,b,c\}$, уфх $a \in R_{mn}(p)$, $b \in R_{m(n-2)}(p)$ ш $c \in R_{m(n-4)}(p)$, хёыш $(m,n,p) \not =(1,8,2),(1,10,2)$; $\{a,b,c\}=\{5,7,17\}$, хёыш $(m,n,p) =(1,8,2)$; $\{a,b,c\}=\{7,17,31\}$, хёыш $(m,n,p) =(1,10,2)$. ┬тшфє ьрыющ ЄхюЁхь√ ╘хЁьр ўшёыр $a$, $b$ ш $c$ эх фхы Є $|Out(S)|$. ╚ч \cite[Єрсы.~2]{Maslova_VasVd2011} ёыхфєхЄ, ўЄю ўшёыр $\{a,b,c\}$ юсЁрчє■Є $3$-ъюъышъє т $\Gamma(S)$, яю¤Єюьє юэш юсЁрчє■Є $3$-ъюъышъє т $\Gamma(G)$ ттшфє ыхьь√ \ref{Maslova_nonadj}. ╦хььр фюърчрэр.

\begin{lem}\label{POmega} ╧єёЄ№ $S \cong P\Omega_{n+1}(q)$, уфх $n \ge 6$ ўхЄэю ш $q=p^m$ эхўхЄэю. ╥юуфр
 $\Gamma(G)$ ёюфхЁцшЄ $3$-ъюъышъє.
\end{lem}

\proof  ┬тшфє \cite[яЁхфыюцхэш ~3.1,~4.3]{Maslova_VasVd2005} ш  \cite[яЁхфыюцхэшх~2.4]{Maslova_VasVd2011} шьххь $\Gamma(PSp_n(q))=\Gamma(P\Omega_{n+1}(q))$. ┬тшфє \cite{Atlas} шьххь $|Out(PSp_n(q))|=|Out(P\Omega_{n+1}(q))|$. ╚ч фюърчрЄхы№ёЄтр ыхьь√ \ref{PSp} ёыхфєхЄ, ўЄю $\Gamma(PSp_n(q))$ яЁш $n \ge 6$ ш эхўхЄэюь $q$ ёюфхЁцшЄ $3$-ъюъышъє $\{a,b,c\}$ Єръє■, ўЄю ўшёыр $a$, $b$ ш $c$ эх фхы Є $|Out(PSp_n(q))|$. ╧Ёшьхэхэшх ыхьь√ \ref{Maslova_nonadj} чртхЁ°рхЄ фюърчрЄхы№ёЄтю.

\begin{lem}\label{POmegapm} ╧єёЄ№ $S \cong P\Omega_n^\varepsilon(q)$, уфх $n \ge 8$ ўхЄэю, $\varepsilon \in \{+,-\}$ ш $q=p^m \ge 2$. ╥юуфр

$(1)$ $\Gamma(G)$ эх ёюфхЁцшЄ $3$-ъюъышъ Єюуфр ш Єюы№ъю Єюуфр, ъюуфр $S \cong P\Omega_8^+(2)${\rm;}

$(2)$ хёыш $S \cong P\Omega_8^+(2)$, Єю ёє∙хёЄтєхЄ ЁрчЁх°шьр  уЁєяяр $H$ Єрър , ўЄю $\Gamma(G)=\Gamma(H)$.
\end{lem}

\proof ┬тшфє \cite[Єрсы.~5.1.A]{Maslova_KlLi} шьххь $|Out(S)|=2\lambda(q^{n/2}-\varepsilon1,4)m$, уфх $\lambda=3$, хёыш $(n,\varepsilon)=(8,+)$ ш $\lambda = 1$ шэрўх. ┼ёыш $(n,q,\varepsilon)=(8,2,+)$, Єю ттшфє \cite[ЄхюЁхьр~1]{ZinMaz} $\Gamma(S)$ ёюёЄюшЄ шч фтєї эхяхЁхёхър■∙шїё  ъышъ, ёыхфютрЄхы№эю ттшфє ыхьь√ \ref{solvable} ёє∙хёЄтєхЄ ЁрчЁх°шьр  уЁєяяр $H$ Єрър , ўЄю $\Gamma(G)=\Gamma(H)$.

╧єёЄ№ $(n,q,\varepsilon)\not =(8,2,+)$. ┬тшфє \cite[Єрсы.~3.5.E,\,3.5.F]{Maslova_KlLi} т $S$ ёє∙хёЄтєхЄ ьръёшьры№эр  яюфуЁєяяр $T$, шчюьюЁЇэр  $P\Omega_{n-1}(q)$. ┬тшфє \cite{Atlas} шьххь $|Out(S)|=|Out(T)|$, хёыш $(n,\varepsilon)\not = (8,+)$, ш $\pi(Out(S)) \subseteq \pi(T) \cup \{3\}$, хёыш $(n,\varepsilon) = (8,+)$. ╩Ёюьх Єюую, хёыш $q$ ўхЄэю, Єю ттшфє \cite{Atlas} шьххь $PSp_n(q) \cong P\Omega_{n+1}(q)$. ╚ч фюърчрЄхы№ёЄт ыхьь \ref{PSp} ш \ref{POmega} ёыхфєхЄ, ўЄю $\Gamma(T)$ ёюфхЁцшЄ $3$-ъюъышъє $\{a,b,c\}$ Єръє■, ўЄю ўшёыр $a$, $b$ ш $c$ эх фхы Є $|Out(T)|$. ╩Ёюьх Єюую, т ёыєўрх $T\cong PSp_6(q)$, уфх $q \not =2$, ўшёыр $a$, $b$ ш $c$ юЄышўэ√ юЄ $3$. ╧Ёшьхэхэшх ыхьь√ \ref{Maslova_nonadj} чртхЁ°рхЄ фюърчрЄхы№ёЄтю.

\begin{lem}\label{Exeptiona} ╧єёЄ№ $S$ --- яЁюёЄр  шёъы■ўшЄхы№эр  уЁєяяр ышхтр Єшяр эрф яюыхь яюЁ фър $q=p^m$. ╥юуфр

$(1)$ $\Gamma(G)$ эх ёюфхЁцшЄ $3$-ъюъышъ Єюуфр ш Єюы№ъю Єюуфр, ъюуфр $S \cong {^3}D_4(2)${\rm;}

$(2)$ хёыш $S \cong {^3}D_4(2)$, Єю ёє∙хёЄтєхЄ ЁрчЁх°шьр  уЁєяяр $H$ Єрър , ўЄю $\Gamma(G)=\Gamma(H)$.
\end{lem}

\proof ╧єёЄ№ $S \cong E_8(q)$. ┬тшфє \cite[Єрсы.~5.1.B]{Maslova_KlLi} шьххь $|Out(S)|=m$. ╚ч \cite[Єрсы.~4]{Maslova_VasVd2011} ёыхфєхЄ, ўЄю ўшёыр $r_{30} \in R_{30m}(p)$, $r_{24} \in R_{24m}(p)$ ш $r_{20} \in R_{20m}(p)$ юсЁрчє■Є $3$-ъюъышъє т $\Gamma(S)$. ┬тшфє ьрыющ ЄхюЁхь√ ╘хЁьр ўшёыр $r_{30}$, $r_{24}$ ш $r_{20}$ эх фхы Є $|Out(S)|$, ттшфє ыхьь√ \ref{Maslova_nonadj} юэш юсЁрчє■Є $3$-ъюъышъє т $\Gamma(G)$.

╧єёЄ№ $S \cong E_7(q)$. ┬тшфє \cite[Єрсы.~5.1.B]{Maslova_KlLi} шьххь $|Out(S)|=(2,q-1)m$. ╚ч \cite[Єрсы.~4]{Maslova_VasVd2011} ёыхфєхЄ, ўЄю ўшёыр $r_{18} \in R_{18m}(p)$, $r_{14} \in R_{14m}(p)$ ш $r_{12} \in R_{12m}(p)$ юсЁрчє■Є $3$-ъюъышъє т $\Gamma(S)$. ┬тшфє ьрыющ ЄхюЁхь√ ╘хЁьр ўшёыр $r_{18}$, $r_{14}$ ш $r_{12}$ эх фхы Є $|Out(S)|$, ттшфє ыхьь√ \ref{Maslova_nonadj} юэш юсЁрчє■Є $3$-ъюъышъє т $\Gamma(G)$.

╧єёЄ№ $S \cong E_6(q)$. ┬тшфє \cite[Єрсы.~5.1.B]{Maslova_KlLi} шьххь $|Out(S)|=2(3,q-1)m$. ╚ч \cite[Єрсы.~4]{Maslova_VasVd2011} ёыхфєхЄ, ўЄю ўшёыр $r_{9} \in R_{9m}(p)$, $r_{8} \in R_{8m}(p)$ ш $r_{5} \in R_{5m}(p)$ юсЁрчє■Є $3$-ъюъышъє т $\Gamma(S)$. ┬тшфє ьрыющ ЄхюЁхь√ ╘хЁьр ўшёыр $r_{9}$, $r_{8}$ ш $r_{5}$ эх фхы Є $|Out(S)|$, яю¤Єюьє ттшфє ыхьь√ \ref{Maslova_nonadj} юэш юсЁрчє■Є $3$-ъюъышъє т $\Gamma(G)$.

╧єёЄ№ $S \cong {^2}E_6(q)$. ┬тшфє \cite[Єрсы.~5.1.B]{Maslova_KlLi} шьххь $|Out(S)|=2(3,q+1)m$. ╚ч \cite[Єрсы.~4]{Maslova_VasVd2011} ёыхфєхЄ, ўЄю ўшёыр шч $r_{18} \in R_{18m}(p)$, $r_{12} \in R_{12m}(p)$ ш $r_{10} \in R_{10m}(p)$ юсЁрчє■Є $3$-ъюъышъє т $\Gamma(S)$. ┬тшфє ьрыющ ЄхюЁхь√ ╘хЁьр ўшёыр $r_{18}$, $r_{12}$ ш $r_{10}$ эх фхы Є $|Out(S)|$, яю¤Єюьє ттшфє ыхьь√ \ref{Maslova_nonadj} юэш юсЁрчє■Є $3$-ъюъышъє т $\Gamma(G)$.

╧єёЄ№ $S \cong {^3}D_4(q)$. ┬тшфє \cite[Єрсы.~5.1.B]{Maslova_KlLi} шьххь $|Out(S)|=3m$. ┼ёыш $q=2$, Єю Єю ттшфє \cite{Atlas} шьххь $\pi(S)=\pi(Aut(S))=\{2,3,7\} \cup \{13\}$ ш ьэюцхёЄтр $\{2,3,7\}$ ш $\{13\}$  ты ■Єё  ъышърьш т $\Gamma(S)$. ╥ръшь юсЁрчюь, ттшфє ыхьь√ \ref{solvable} ёє∙хёЄтєхЄ ЁрчЁх°шьр  уЁєяяр $H$ Єрър , ўЄю $\Gamma(G)=\Gamma(H)$. ┼ёыш $q=4$, Єю ўшёыр $7 \in R_{3}(4)$, $13 \in R_{6}(4)$ ш $241 \in R_{12}(4)$ юсЁрчє■Є $3$-ъюъышъє т $\Gamma(S)$ ш эх фхы Є $|Out(S)|=6$, яю¤Єюьє ттшфє ыхьь√ \ref{Maslova_nonadj} юэш юсЁрчє■Є $3$-ъюъышъє т $\Gamma(G)$.  ╧Ёхфяюыюцшь, ўЄю $q>2$ ш $q \not = 4$. ╥юуфр шч \cite[Єрсы.~4]{Maslova_VasVd2011} ёыхфєхЄ, ўЄю ўшёыр шч $r_3 \in R_{3m}(p)$, $r_6 \in R_{6m}(p)$ ш $r_{12} \in R_{12m}(p)$ юсЁрчє■Є $3$-ъюъышъє т $\Gamma(S)$. ┬тшфє ьрыющ ЄхюЁхь√ ╘хЁьр ўшёыр $r_{3}$, $r_{6}$ ш $r_{12}$ эх фхы Є $|Out(S)|$, яю¤Єюьє ттшфє ыхьь√ \ref{Maslova_nonadj} юэш юсЁрчє■Є $3$-ъюъышъє т $\Gamma(G)$.

╧єёЄ№ $S \cong F_4(q)$. ┬тшфє \cite[Єрсы.~5.1.B]{Maslova_KlLi} шьххь $|Out(S)|=(2,p)m$. ╚ч \cite[Єрсы.~4]{Maslova_VasVd2011} ёыхфєхЄ, ўЄю ўшёыр $r_{12} \in R_{12m}(p)$, $r_{8} \in R_{8m}(p)$ ш $r_{4} \in R_{4m}(p)$ юсЁрчє■Є $3$-ъюъышъє т $\Gamma(S)$. ┬тшфє ьрыющ ЄхюЁхь√ ╘хЁьр ўшёыр $r_{12}$, $r_{8}$ ш $r_{4}$ эх фхы Є $|Out(S)|$, ттшфє ыхьь√ \ref{Maslova_nonadj} юэш юсЁрчє■Є $3$-ъюъышъє т $\Gamma(G)$.

╧єёЄ№ $S \cong {^2}F_4(2^{2n+1})$, уфх $n \ge 1$. ┬тшфє \cite[Єрсы.~5.1.B]{Maslova_KlLi} шьххь $|Out(S)|=2n+1$. ╧єёЄ№ $n \ge 2$. ╚ч \cite[Єрсы.~4]{Maslova_VasVd2011} ёыхфєхЄ, ўЄю ы■с√х ўшёыр $s_2 \in R_{4n+2}(2) \subset \pi(2^{2n+1}+1) \setminus \{3\} $, $s_{3} \in R_{8n+4}(2) \subset \pi(2^{4n+2}+1) \setminus \{3\}$ ш $s_{4} \in R_{12n+6}(2) \subset \pi(2^{4n+2}-2^{2n+1}+1) \setminus \{3\}$ юсЁрчє■Є $3$-ъюъышъє т $\Gamma(S)$. ┬тшфє ьрыющ ЄхюЁхь√ ╘хЁьр ўшёыр $s_2$, $s_3$ ш $s_4$ эх фхы Є $2n+1$, ЄхяхЁ№ шч ыхьь√ \ref{Maslova_nonadj} ёыхфєхЄ, ўЄю юэш юсЁрчє■Є $3$-ъюъышъє т $\Gamma(G)$.
╧єёЄ№ $n=1$. ╚ч \cite[Єрсы.~4]{Maslova_VasVd2011} ёыхфєхЄ, ўЄю ўшёыр $2$, $19$ ш $37$ юсЁрчє■Є $3$-ъюъышъє т $\Gamma(S)$. ╧юёъюы№ъє ¤Єш ўшёыр эх фхы Є $|Out(S)|$, шч ыхьь√ \ref{Maslova_nonadj} ёыхфєхЄ, ўЄю юэш юсЁрчє■Є $3$-ъюъышъє т $\Gamma(G)$.

╧єёЄ№ $S \cong {^2}F_4(2)'$. ┬тшфє \cite{Atlas} шьххь $|Out(S)|=2$. ╚ч \cite[Єрсы.~4]{Maslova_VasVd2011} ёыхфєхЄ, ўЄю ўшёыр $3$, $5$ ш $13$ юсЁрчє■Є $3$-ъюъышъє т $\Gamma(S)$. ╧юёъюы№ъє ¤Єш ўшёыр эх фхы Є $|Out(S)|$, шч ыхьь√ \ref{Maslova_nonadj} ёыхфєхЄ, ўЄю юэш юсЁрчє■Є $3$-ъюъышъє т $\Gamma(G)$.

╧єёЄ№ $S \cong G_2(q)$. ┬тшфє \cite[Єрсы.~5.1.B]{Maslova_KlLi} шьххь $|Out(S)|=2^{\frac{(3,p)-1}{2}}m$. ╧єёЄ№ $r_{6} \in R_{6m}(p)$, $r_{3} \in R_{3m}(p)$ ш $r_{2} \in R_{2m}(p)$, хёыш $(m,p)\not = (2,2),(3,2)$, $r_6=13$, $r_3=7$ ш $r_2=5$, хёыш $(m,p)=(2,2)$, $r_6=19$, $r_3=73$ ш $r_2=7$, хёыш $(m,p)=(3,2)$. ╚ч \cite[Єрсы.~4]{Maslova_VasVd2011} ёыхфєхЄ, ўЄю ўшёыр $\{r_{6},r_{3},r_{2}\}$ юсЁрчє■Є $3$-ъюъышъє т $\Gamma(S)$. ┬тшфє ьрыющ ЄхюЁхь√ ╘хЁьр ўшёыр $r_{6}$, $r_{3}$ ш $r_{2}$ эх фхы Є $|Out(S)|$, яю¤Єюьє ттшфє ыхьь√ \ref{Maslova_nonadj} юэш юсЁрчє■Є $3$-ъюъышъє т $\Gamma(G)$.

╧єёЄ№ $S \cong {^2}G_2(3^{2n+1})$, уфх $n \ge 1$. ┬тшфє \cite[Єрсы.~5.1.B]{Maslova_KlLi} шьххь $|Out(S)|=2n+1$. ╚ч \cite[Єрсы.~4]{Maslova_VasVd2011} ёыхфєхЄ, ўЄю ы■с√х ўшёыр $s_1 \in R_{2n+1}(3) \subset \pi(3^{2n+1}-1)$, $s_2 \in R_{4n+2}(3) \subset \pi(3^{2n+1}+1)$ ш $s_{3} \in R_{12n+6}(3) \subset \pi(3^{4n+2}-3^{2n+1}+1) $ юсЁрчє■Є $3$-ъюъышъє т $\Gamma(S)$. ┬тшфє ьрыющ ЄхюЁхь√ ╘хЁьр ўшёыр $s_1$, $s_2$ ш $s_3$ эх фхы Є $2n+1$, ЄхяхЁ№ шч ыхьь√ \ref{Maslova_nonadj} ёыхфєхЄ, ўЄю юэш юсЁрчє■Є $3$-ъюъышъє т $\Gamma(G)$.

╧єёЄ№ $S \cong {^2}B_2(2^{2n+1})$, уфх $n \ge 1$. ┬тшфє \cite[Єрсы.~5.1.B]{Maslova_KlLi} шьххь $|Out(S)|=2n+1$. ╚ч \cite[Єрсы.~4]{Maslova_VasVd2011} ёыхфєхЄ, ўЄю ўшёыр $2$, $s_1 \in R_{2n+1}(2) \subset \pi(2^{2n+1}-1)$ ш $s_2 \in R_{8n+4}(2) \subset \pi(2^{4n+2}+1)$ юсЁрчє■Є $3$-ъюъышъє т $\Gamma(S)$. ┬тшфє ьрыющ ЄхюЁхь√ ╘хЁьр ўшёыр $s_1$ ш $s_2$ эх фхы Є $2n+1$, ЄхяхЁ№ шч ыхьь√ \ref{Maslova_nonadj} ёыхфєхЄ, ўЄю ўшёыр $\{2, s_1, s_2\}$ юсЁрчє■Є $3$-ъюъышъє т $\Gamma(G)$. ╦хььр фюърчрэр.

\medskip ─юърчрЄхы№ёЄтю ЄхюЁхь \ref{ASasS} ш \ref{AlmostSimple} ёыхфєхЄ шч ыхьь \ref{Sporadic}--\ref{Exeptiona}.

\bigskip

╤тхфхэш  юс ртЄюЁрї:

\smallskip

├╬╨╪╩╬┬ ╚ы№  ┴юЁшёютшў,

\smallskip

╨╬╤╤╚▀, у. ┼ърЄхЁшэсєЁу, єы. ╤.╩ютрыхтёъющ, 16, 620990,

╚эёЄшЄєЄ ьрЄхьрЄшъш ш ьхїрэшъш шь. ═.═. ╩Ёрёютёъюую ╙Ё╬ ╨└═,

х-mail: ilygor8@gmail.com

\smallskip

╠└╤╦╬┬└ ═рЄры№  ┬ырфшьшЁютэр,

\smallskip

╨╬╤╤╚▀, у. ┼ърЄхЁшэсєЁу, єы. ╤.╩ютрыхтёъющ, 16, 620990,

╚эёЄшЄєЄ ьрЄхьрЄшъш ш ьхїрэшъш шь. ═.═. ╩Ёрёютёъюую ╙Ё╬ ╨└═,

╨╬╤╤╚▀, у. ┼ърЄхЁшэсєЁу, єы. ╠шЁр, 19, 620002,

╙Ёры№ёъшщ ЇхфхЁры№э√щ єэштхЁёшЄхЄ,

х-mail: butterson@mail.ru

\end{document}